\documentclass[12pt,onefignum,onetabnum,a4paper]{article}

\usepackage{graphicx}
\usepackage{amsmath,amsfonts,amssymb}  
\usepackage{url}  
\usepackage{psfrag}
\usepackage{color,caption}

\newcommand{\mR}{\mathbb{R}}
\newcommand{\mT}{\mathbb{T}}

\newcommand{\fref}[1]{Fig.~\ref{#1}}

\newcommand{\Fref}[1]{Figure~\ref{#1}}

\newcommand{\sref}[1]{Section~\ref{#1}}
\newcommand{\Sref}[1]{Section~\ref{#1}}

\newcommand{\xeq}{\mathbf{x}^*}
\newcommand{\iso}[1]{I(\gamma_{#1})}
\newcommand{\bu}{\mathbf{u}}
\newcommand{\bv}{\mathbf{v}}
\newcommand{\bw}{\mathbf{w}}

\newcommand{\bx}{\mathbf{x}}

\newcommand{\sPO}[1]{$\Gamma_p$}

\begin{document}



\title{A continuation approach to \\ computing phase resetting curves}

\author{Peter Langfield\footnote{{\tt peter.langfield@u-bordeaux.fr}},
Bernd Krauskopf$\,^{\dagger}$ and Hinke M. Osinga\footnote{{\tt b.krauskopf@auckland.ac.nz}; {\tt h.m.osinga@auckland.ac.nz}}\\[4mm]
$^*$IHU Liryc, Electrophysiology and Heart Modeling Institute, \\ Fondation Bordeaux Universit{\'e}, F-33600, Pessac- Bordeaux, France\\[2mm]
$^{\dagger}$ Department of Mathematics, The University of Auckland, \\ Private Bag 92019, Auckland 1142, New Zealand}

\date{March 2020}

\maketitle

%
%
%


\abstract{Phase resetting is a common experimental approach to investigating the behaviour of oscillating neurons. Assuming repeated spiking or bursting, a phase reset amounts to a brief perturbation that causes a shift in the phase of this periodic motion. The observed effects not only depend on the strength of the perturbation, but also on the phase at which it is applied. The relationship between the change in phase after the perturbation and the unperturbed old phase, the so-called phase resetting curve, provides information about the type of neuronal behaviour, although not all effects of the nature of the perturbation are well understood. In this chapter, we present a numerical method based on the continuation of a multi-segment boundary value problem that computes phase resetting curves in ODE models. Our method is able to deal effectively with phase sensitivity of a system, meaning that it is able to handle extreme variations in the phase resetting curve, including resets that are seemingly discontinuous. We illustrate the algorithm with two examples of planar systems, where we also demonstrate how qualitative changes of a phase resetting curve can be characterised and understood. A seven-dimensional example emphasises that our method is not restricted to planar systems, and illustrates how we can also deal with non-instantaneous, time-varying perturbations.}

%
\section{Introduction}
Measuring phase resetting is a common approach for testing neuronal responses in experiments: a brief current injection perturbs the regular spiking behaviour of a neuron, resulting generally in a shifted phase as the neuron returns to its regular oscillating behaviour. This phase shift can be advanced or delayed---meaning that the next spike arrives earlier or later compared with the unperturbed spiking oscillation---and which effect occurs also depends on the moment when the current is applied; see~\cite{ErmentroutTerman2010} for more details. A plot of the shifted phase $\vartheta_{\rm new}$ versus the original phase $\vartheta_{\rm old}$ at which the current was applied is known as the \emph{phase transition curve} (PTC). Experimentally, it is often easier to represent the reset in terms of the resulting phase difference $\vartheta_{\rm new} - \vartheta_{\rm old}$ as a function of $\vartheta_{\rm old}$, which can be measured as the time to the next spike; such a representation is called a \emph{phase response curve} or \emph{phase resetting curve} (PRC). 

The shape of a PTC or PRC of a given system obviously depends on the size of the applied perturbation: already for quite small amplitudes, nonlinear effects can dramatically affect a PTC or PRC. The shape of the PTC or PRC has been used to classify neuronal behaviour~\cite{brown2004, ermentrout96, hansel1995}, where the underlying assumption is that the size of the applied perturbation is sufficiently small. Hodgkin~\cite{hodgkin1948} distinguished between so-called Type-I and Type-II excitable membranes, where neurons with membranes of Type II are not able to fire at arbitrarily low frequencies. Note that transitions from Type-I to Type-II can occur when system parameters are changed~\cite{ermentrout12}. Ermentrout~\cite{ermentrout96} found that the PRC of a Type-I neuron always has the same sign, while that of a Type-II neuron changes sign; this means that the PTC is always entirely above or below the diagonal for Type-I neurons, while it intersects the diagonal for Type-II neurons. In either case, the PTC is invertible for sufficiently small perturbation amplitudes, since it can be viewed as a continuous and smooth deformation of the identity, which is the PTC in the limit of zero amplitude. Invertibility itself has also been used as a distinguishing property of PTCs: noninvertible PTCs are said to be of type~0 (or strong) and invertible PTCs are of type~1 (or weak)~\cite{ErmentroutTerman2010, Winfree2001}. If an increasingly stronger perturbation is applied, for example, in the context of synchronisation, it is well known that PTCs can change from type~1 to type~0, that is, become noninvertible~\cite{glass84, Winfree2001}. 

A motivation in recent work on phase resetting has been the idea of interpreting the PTC as defining a one-dimensional phase-reduction model that, hopefully, captures the essential dynamics of a possibly high-dimensional oscillating system. The main interest is in coupled systems, formed by two or more (planar) systems with known PRCs; for example, see~\cite{PietrasDaffertshofer2019, Schultheiss2012} for mathematical as well as experimental perspectives. Unfortunately, the convergence back to the limit cycle after some perturbation can be quite slow for a coupled system, such that only (infinitesimally) small perturbations are accurately described. Furthermore, it makes physiological sense to assume a time-varying input, usually in the form of a short input pulse, rather than the instantaneous perturbation assumed for the theoretical phase reset. Moreover, the perturbation may be repeated at regular intervals. In this context, PTCs and PRCs can be useful for explaining the resetting behaviour, though strictly speaking, the theory is only valid at low firing rates~\cite{ermentrout05, WilsonErmentrout2019}. More recently, the idea of a phase-amplitude description has led to a better understanding of the effects resulting from these kinds of repeated time-varying resets~\cite{CastejonGuillamon2020, castejon13, MongaWilsonMatchenMoehlis2019, PerezSearaHuguet2019, Wedgwood2013}.

From a dynamical systems perspective, the key question of phase resetting is how the perturbed initial conditions relax back to an attracting periodic orbit $\Gamma$ with period $T_\Gamma$ of an underlying continuous-time model, which we take here to be a vector field on $\mR^n$, that is, a system of $n$ first-order autonomous ordinary differential equations. All points in its basin $\mathcal{B}(\Gamma)$ converge to $\Gamma$, and they do so with a given asymptotic phase. The subset of all points in $\mathcal{B}(\Gamma)$ that converge to $\Gamma$ in phase with the point $\gamma_{\vartheta} \in \Gamma$, where $\vartheta \in [0, 1)$ by convention, is called the (forward-time) \emph{isochron} of $\gamma_{\vartheta}$, which we refer to as $\iso{\vartheta}$. Isochrons were defined and named by Winfree~\cite{Winfree1974}. Guckenheimer~\cite{guckenheimer75} showed that $\iso{\vartheta}$ is, in fact, an $(n-1)$-dimensional invariant stable manifold of the attracting fixed point $\gamma_{\vartheta} \in \Gamma$ under the time-$T_\Gamma$ map. In particular, it follows that $\iso{\vartheta}$ is tangent to the attracting linear eigenspace of $\gamma_{\vartheta}$ and, hence, transverse to $\Gamma$.  Moreover, the $\vartheta$-dependent family of all isochrons $\iso{\vartheta}$ foliates the basin ${\mathcal B}(\Gamma)$. In other words, any point in ${\mathcal B}(\Gamma)$ has a unique asymptotic phase determined by the isochron it lies on.

For a given $\vartheta_{\rm old}$, consider now the perturbed point $\gamma_{\vartheta_{\rm old}} + A \, \mathbf{d} \in {\mathcal B}(\Gamma)$, obtained from $\gamma_{\vartheta_{\rm old}} \in \Gamma$ by applying the perturbation of strength $A$ in the given direction $\mathbf{d}$. The asymptotic phase $\vartheta_{\rm new}$ is, hence, uniquely determined by the isochron $\iso{\vartheta_{\rm new}}$ on which this point lies. This defines a circle map $P : [0, 1) \to [0, 1)$ with $P(\vartheta_{\rm old}) = \vartheta_{\rm new}$. Therefore, finding the PTC is equivalent to determining how the perturbed cycle $\Gamma + A \, \mathbf{d} = \{ \gamma_{\vartheta_{\rm old}} + A \, \mathbf{d} \;|\; \vartheta_{\rm old} \in [0, 1) \}$ intersects the foliation of ${\mathcal B}(\Gamma)$ by the isochrons $\iso{\vartheta_{\rm new}}$ for $\vartheta_{\rm new} \in [0, 1)$. Notice further that the PTC is the graph of the circle map $P$ on the unit torus $\mathbb{T}^2$, represented by the unit square $[0, 1) \times [0, 1)$.

When considering the amplitude A of the perturbation as a parameter (while keeping the direction $\mathbf{d}$ fixed throughout), one can deduce some important properties of the associated PTC. Suppose that $0 < A_{\rm max}$ is such that $\Gamma_A := \Gamma + A \, \mathbf{d} \in {\mathcal B}(\Gamma)$ for all $0 \leq A \leq A_{\rm max}$. Then none of these perturbed cycles $\Gamma_A$ intersects the boundary of the basin ${\mathcal B}(\Gamma)$ and the associated circle map $P = P_A$ is well defined for all $\vartheta_{\rm old} \in [0, 1)$. The map $P_0$ for zero perturbation amplitude is the identity on $\mathbb{T}^2$, which means that, as its graph, the PTC is the diagonal on $[0, 1) \times [0, 1)$  and a 1:1 torus knot on $\mathbb{T}^2$; in particular, $P_0$ is invertible, that is, it is injective and surjective. Because of smooth dependence on the amplitude $A$ and the fact that $P_A$ is a function over $[0,1)$, the PTC remains a 1:1 torus knot on $\mT^2$ and $P_A$ is surjective for all $0 \leq A \leq A_{\rm max}$. 

Since the isochrons are transverse to $\Gamma$, the circle map $P_A$ is $C^1$-close to the identity, and hence, also injective, for sufficiently small $A$. As the graph of a near-identity transformation, the PTC is then strictly monotone, invertible, and hence, of type~1 (or weak) in the notation of~\cite{ErmentroutTerman2010, Winfree2001}. While surjectivity is preserved, injectivity may be lost before $A = A_{\rm max}$ is reached. Indeed, the PTC is either invertible for all $0 \leq A \leq A_{\rm max}$, or there is a maximal $0 < A_{\rm inv} < A_{\rm max}$ such that $P_A$ is invertible only for all $0 < A \leq A_{\rm inv}$. The loss of injectivity of $P_A$ at $A = A_{\rm inv}$ happens generically because of the emergence of an inflection point. For $A_{\rm inv} < A \leq A_{\rm max}$ this transition creates a local minimum and a local maximum of the PTC, which is now no longer invertible and so of type~0 (or strong)  in the notation of~\cite{ErmentroutTerman2010, Winfree2001}. As we will show, an inflection point of $P_A$ corresponds to a cubic tangency between the perturbed cycle $\Gamma_A$ and an isochron. Indeed, additional inflection points and, hence, local minima and maxima may appear at subsequent cubic isochron tangencies. Since $P_A$ is a circle map, these must come in pairs; hence, counting the number of its local maxima (or minima) would provide a further refinement of the notation of a type~0 (or strong) PTC.

The above discussion shows that, when the applied perturbation $A$ is sufficiently weak, it suffices to consider only the linear approximation to the isochron family, which is given by the $\vartheta$-family of stable eigenspaces of the time-$T_\Gamma$ map for each $\vartheta$. In practice, nonlinear effects are essential, especially when multiple time scales are present or the phase reset involves relatively strong perturbations. Isochrons are often highly nonlinear objects of possibly very complicated geometry~\cite{lko-chaos2014, Winfree1974}. While the geometric idea of isochrons determining the phase resets has been around since the mid 1970s, the practical implementation has proven rather elusive. In practice, it is not at all straightforward to compute the isochrons of a periodic orbit. In planar systems, when such isochrons are curves, three different approaches have been proposed, based on Fourier averages~\cite{Mauroy2012, mauroy14}, a parametrisation formulated in terms of a functional equation~\cite{guillamon2009, huguet2013}, and continuation of solutions to a suitably posed two-point boundary value problem~\cite{lko-chaos2014, om-siads2010}. In principle, all three approaches generalise to higher-dimensional isochrons, but there are only few explicit examples~\cite{guillamon2009, mauroy14}.

From the knowledge of the isochron foliation of ${\cal B}(\Gamma)$, one can immediately deduce geometrically the phase resetting for perturbations of any strength and in any direction. However, already for planar and certainly for higher-dimensional systems, this is effectively too much information when one is after the PTC resulting from a perturbation in a fixed direction and with a specific amplitude. In essence, finding a PTC or PRC remains the one-dimensional problem of finding the asymptotic phase of all points on the perturbed cycle. 

In this chapter, we show how this can be achieved with a multi-segment boundary value problem formulation. Specifically, we adapt the approach from~\cite{lko-chaos2014, lko-siads2015} to set up the calculation of the circle map $P_A$ by continuation, first in $A$ from $A=0$ for fixed $\vartheta_{\rm old}$, and then in $\vartheta_{\rm old} \in [0, 1)$ for fixed $A$. In this way, we obtain accurate numerical approximations of the PTC or PRC as continuous curves, even when the system shows strong phase sensitivity. The set-up is extremely versatile, and the direct computation of a PTC in this way does not require the system to be planar. We demonstrate our method with a constructed example going back to Winfree~\cite[Chapter~6]{Winfree2001}, where we also show how injectivity is lost in a first cubic tangency of $\Gamma_A$ with an isochron. The robustness of the method is then illustrated with the computation of a PTC of a perturbed cycle that cuts through a region of extreme phase sensitivity in the (planar) FithHugh--Nagumo system; in spite of very large derivatives due to this phase sensitivity, the PTC is computed accurately as a continuous curve. Our final example of a seven-dimension system from~\cite{TKM}  modelling a type of cardiac pacemaker cell shows that our approach also works in higher dimensions; this system also features phase sensitivity due to the existence of different time scales.

This chapter is organised as follows. In the next section, we provide precise details of the setting and explain the definitions used. \Sref{sec:PRC} presents the numerical set-up for computing a resetting curve by continuation of a multi-segment boundary value problem. We then discuss two planar examples in depth, which are both taken from~\cite{lko-siads2015}: a variation of Winfree's model in \sref{sec:winfree} and the FitzHugh--Nagumo system in \sref{sec:fhn}. The third and higher-dimensional example from~\cite{TKM} is presented in \sref{sec:TKM}. A summary of the results is given in \sref{sec:conclusions}, where we also discuss some consequences of our findings and directions of future research.

%
\section{Basic setting and definitions}
\label{sec:settings}
As mentioned in the introduction, we consider a dynamical system with an attracting periodic orbit $\Gamma$. For simplicity, we assume that the state space is $\mR^n$ and consider the dynamical system
\begin{equation}
\label{eq:basic}
  \dot{\bx} = \mathbf{F}(\bx),
\end{equation}
where $\mathbf{F} : \mR^n \to \mR^n$ is at least once continuously differentiable. We assume that system~\eqref{eq:basic} has an attracting periodic orbit $\Gamma$ with period $T_\Gamma$, that is,
\begin{displaymath}
  \Gamma := \{ \gamma(t) \in \mR^n \;|\; 0 \leq t \leq T_\Gamma \mbox{ with } \gamma(T_\Gamma) = \gamma(0) \},
\end{displaymath}
and $T_\Gamma$ is minimal with this property. We associate a phase $\vartheta \in [0, 1)$ with each point $\gamma_\vartheta \in \Gamma$, defining $\gamma_\vartheta := \gamma(t)$ with $t = \vartheta T_\Gamma$. Here $\gamma_0 := \gamma(0)$ needs to be chosen, which is usually done by fixing it to correspond to a maximum in the first component. The (forward-time) isochron $\iso{\vartheta}$ associated with $\gamma_\vartheta \in \Gamma$ is then defined in terms of initial conditions $\bx(0)$ of forward trajectories $\bx := \{ \bx(t) \in \mR^n \;|\; t \in \mR \}$ of system~\eqref{eq:basic} that accumulate on $\Gamma$, namely, as
\begin{displaymath}
  \iso{\vartheta} := \{ \bx(0) \in \mR^n \;|\; \lim_{s \to \infty} \bx(s \, T_\Gamma) = \gamma_\vartheta \}.
\end{displaymath}
In other words, the trajectory $\bx$ approaches $\Gamma$ in phase with $\gamma_\vartheta$. Note that $\iso{\vartheta}$ is the stable manifold of the fixed point $\gamma_\vartheta$ of the time-$T_\Gamma$ return map; in particular, this means that $\iso{\vartheta}$ is of dimension $n-1$ and tangent at $\gamma_\vartheta$ to the stable eigenspace $E(\gamma_\theta)$, which is part of the stable Floquet bundle of $\Gamma$~\cite{guckenheimer75}; we utilise this property when computing isochrons, and also when computing a PTC or PRC.

We are now ready to give formal definitions of the PTC and PRC; see also~\cite{ErmentroutTerman2010}.
\textbf{Definition}[Phase Transition Curve]
\label{def:PTC}
The phase transition curve or PTC associated with a perturbation of amplitude $A \geq 0$ in the direction $\mathbf{d} \in \mR^n$ is the graph of the map $P : [0, 1) \to [0, 1)$ defined as follows. For $\vartheta \in [0, 1)$, the image $P(\vartheta)$  is the phase $\varphi$ associated with the isochron $\iso{\varphi}$ that contains the point $\gamma_\vartheta + A \, \mathbf{d}$ for $\gamma_\vartheta \in \Gamma$.
%
\textbf{Definition}[Phase Response Curve]
\label{def:PRC}
The phase response curve or PRC associated with a perturbation of amplitude $A \geq 0$ in the direction $\mathbf{d} \in \mR^n$ is the graph of the phase difference $\Delta(\vartheta) = P(\vartheta) - \vartheta \, ({\rm mod} \ 1)$, where the map $P$ is as above.
%

The definitions of the PTC and PRC are based on knowledge of the (forward-time) isochron $\iso{\varphi}$ associated with a point $\gamma_\varphi \in \Gamma$. We previously designed an algorithm based on continuation of a two-point boundary value problem (BVP) that computes one-dimensional (forward-time and backward-time) isochrons of a planar system up to arbitrarily large arclengths~\cite{lko-chaos2014, lko-siads2015, om-siads2010}. Here, we briefly describe this algorithm in its simplest form, because this is useful for understanding the basic set-up, and for introducing some notation. The description is presented in the style that is used for implementation in the software package {\sc Auto}~\cite{Doedel1981, auto}. In particular, we consider a time-rescaled version of the vector field~\eqref{eq:basic}, which represents an orbit segment $\{ \bx(t)\ \;|\; 0 \leq t \leq T \}$ of~\eqref{eq:basic} as the orbit segment $\{ \bu(t) \;|\; 0 \leq t \leq 1 \}$ of the vector field
\begin{equation}
\label{eq:rsvf}
  \dot{\bu} = T \, \mathbf{F}(\bu), 
\end{equation}
so that the total integration time $T$ is now a parameter of the system. 

We approximate $\iso{0}$ as the set of initial points of orbit segments that end on the linear space $E(\gamma_0)$, the linearised isochron of $\iso{0}$, close to $\gamma_0$ after integer multiples of the period $T_\Gamma$. These points are formulated as initial points $\bu(0)$ of orbit segments $\bu$ that end on $E(\gamma_0)$ at a distance $\eta$ from $\gamma_0$; hence, $\eta$ defines a one-parameter family of orbit segments. Each orbit segment in this family is a solution of system~\eqref{eq:rsvf} with $T = k \, T_\Gamma$ for $k \in \mathbb{N}$; the corresponding boundary conditions are:
\begin{equation}
\label{eq:onv}
  \left[ \bu(1) - \gamma_0 \right] \cdot \bv_0^\perp = 0,
\end{equation}
\begin{equation}
\label{eq:eta}
  \left[ \bu(1) - \gamma_0 \right] \cdot \bv_0^{\phantom{\perp}} = \eta,
\end{equation}
where $\bv_0$ is the normalised vector that spans $E(\gamma_0)$ and $\bv_0^\perp$ is perpendicular to it. Note that $\Gamma$ itself, when starting from $\gamma_0$, is a solution to the two-point BVP~\eqref{eq:rsvf}--\eqref{eq:eta} with $T = T_\Gamma$ and $\eta = 0$. This gives us a first solution to start the continuation for computing $\iso{0}$. We fix $T = T_\Gamma$ and continue the orbit segment $\bu$ in $\eta$ up to a maximum prespecified tolerance $\eta = \eta_{\rm max}$. As the end point $\bu(1)$ is pushed away from $\gamma_0$ along $E(\gamma_0)$, the initial point $\bu(0)$ traces out a portion of $\iso{0}$. 

Once we reach $\eta = \eta_{\rm max}$, we can extend $\iso{0}$ further by considering points that map to $E(\gamma_0)$ after one additional period, that is, after time $T = 2 \, T_\Gamma$. We start the continuation with the orbit segment formed by concatenation of the final orbit segment with $\Gamma$; here, we rescale time such that this first orbit is again defined for $0 \leq t \leq 1$, we set $T = 2 \, T_\Gamma$, and $\eta = 0$. Note that this orbit segment has a discontinuity at $t = \frac{1}{2}$, but it is very small and {\sc Auto} will automatically correct and close it as part of the first continuation step. This correction will cause a small shift in $\eta$ away from $0$, but $\eta$ will still be much smaller than $\eta_{\rm max}$ (in absolute value). We can keep extending $\iso{0}$ further in this way, by continuation with $T = k \, T_\Gamma$, for integers $k > 2$. See~\cite{lko-chaos2014, om-siads2010} for more details on the implementation and, in particular, see~\cite{krauskopf08, om-siads2010} for details on how to find $E(\gamma_0)$ represented by the first vector $\bv_0$ in the stable Floquet bundle of $\Gamma$.

The computational set-up forms a well-posed two-point BVP with a one-parameter solution family that can be found by continuation, provided the following equality holds for the dimension {\sf NDIM} of the problem, the number {\sf NBC} of boundary conditions, and the number {\sf NPAR}  of free parameters: ${\sf NDIM} - {\sf NBC} + {\sf NPAR} = 1$. Indeed, for the computation of $\iso{0}$, we have ${\sf NDIM} = 2$, because we assumed that the system is planar; ${\sf NBC} = 2$, namely, one condition to restrict $\bu(1)$ to the linearised isochron of $\iso{0}$, and one condition to fix its distance to $\gamma_0$; and ${\sf NPAR} = 1$, because we free the parameter $\eta$.

To compute $\iso{\varphi}$ for other $\varphi \in [0, 1)$, this same approach can be used, working with a shifted periodic orbit $\Gamma$ so that its head point is $\gamma_\varphi$, and determining the associated direction vector $\bv_\varphi$ that spans the eigenspace $E(\gamma_\varphi)$ to which $\iso{\varphi}$ is tangent. In~\cite{om-siads2010}, approximations of $\gamma_\varphi$ and $\bv_\varphi$ are obtained by interpolation of the respective mesh discretisations from {\sc Auto}. We describe an alternative approach in~\cite{lko-chaos2014}, where we consider $\iso{\varphi}$ as the set of initial points of orbit segments that end in the linear space $E(\gamma_0)$ of $\iso{0}$ sufficiently close to $\gamma_0$ after total integration time $T =  k \, T_\Gamma + (1 - \varphi) \, T_\Gamma$.

For the computation of a resetting curve, we use a combination of these two approaches, but rather than interpolation, we shift the periodic orbit by imposing a separate two-point BVP. More precisely, we set up a multi-segment BVP comprised of several subsystems of two-point BVPs; the set-up for this extended BVP is explained in detail in the next section.

%
\section{Algorithm for computing a phase resetting curve}
\label{sec:PRC}
Based on the definition of PTC and PRC, one could now calculate a sufficiently large number of isochrons and determine the resetting curve numerically from data. We prefer to compute the PTC or PRC directly with a BVP set-up and continuation in a very similar way. The major benefit of such a direct approach is that it avoids accuracy restrictions arising from the selection of computed isochrons; in particular, any phase sensitivity of the PTC or PRC will be dealt with automatically as part of the pseudo-arclength continuation with {\sc Auto}~\cite{Doedel1981, auto}. 

For ease of presentation, we will formulate and discuss our continuation set-up for the case of a planar system. We remark, however, that it can readily be extended for use in $\mR^n$ with $n > 2$, because the dimensionality of the problem is not determined by the dimension $n-1$ of the isochrons but by the dimension of the PTC or PRC, which is always one; see also the example in \sref{sec:TKM}. 

The essential difference between calculating a resetting curve rather than an isochron is the following: for an isochron $\iso{\vartheta}$, we compute orbit segments with total integration time $T = T_\Gamma$ (or integer multiples), where we move the end point $\bu(1)$ along the linear approximation of $\iso{\vartheta}$ to some distance $\eta$ from $\Gamma$, while the initial point $\bu(0)$ traces out a new portion of $\iso{\vartheta}$; imagining the same set-up, if we move $\bu(0)$ transverse to $\iso{\vartheta}$, the end point $\bu(1)$ will move to lie on the linearisation of an isochron $\iso{\varphi}$ with a different phase $\varphi$. (Here, one should expect that the distance to $\Gamma$ also changes, but we assume it is still less than $\eta_{\rm max}$).  The key idea behind our approach is that we find a way to determine the different phase $\varphi$, or the phase shift $\varphi - \vartheta$, by allowing $\Gamma$ and its corresponding stable Floquet bundle to rotate as part of an extended system. We ensure the head point of $\Gamma$ moves with the phase-shifted point, that is, the first point on $\Gamma$ will be $\gamma_\varphi$. In this way, we can determine the shifted phase $\varphi$ along any prescribed arc traced out by $\bu(0)$, provided it lies in the basin of attraction of $\Gamma$. For the PTC or PRC associated with a perturbation of amplitude $A \geq 0$ in the direction $\mathbf{d} \in \mR^n$, this arc should be the perturbed cycle $\Gamma + A \, \mathbf{d}$, that is, the closed curve $\{\gamma_\vartheta + A \, \mathbf{d} \;|\; \vartheta \in [0, 1) \}$.

%
\subsection{Continuation set-up for rotated representation of $\Gamma$}
\label{sec:rotate}
We formulate an extended BVP that represents a rotated version of $\Gamma$ with a particular phase, meaning that we automatically determine the phase of the head point relative to $\gamma_0$. To this end, we assume that the zero-phase point $\gamma_0 \in \Gamma$ and its associated linear vector $\bv_0$, or more practical, its perpendicular $\bv_0^\perp$, are readily accessible as stored parameters, or constants that do not change. Hence, even when $\Gamma$ is rotated and its first point is $\gamma_\varphi$ for some different $\varphi \in [0, 1)$, we can still access the coordinates of $\gamma_0$ and $\bv_0^\perp$ from the parameter/constants list.

The extended BVP consists of three components, one to define $\Gamma$, one to define the associated (rotated) linear bundle, and one to define the associated phase. We start by representing $\Gamma$ as a closed orbit segment $\mathbf{g}$ that solves system~\eqref{eq:rsvf} for $T = T_\Gamma$.  Hence, we define
\begin{equation}
\label{eq:po}
  \dot{\mathbf{g}} = T_\Gamma \, \mathbf{F}(\mathbf{g}),
\end{equation}
with periodic boundary condition
\begin{equation}
\label{eq:pobc}
   \mathbf{g}(1) - \mathbf{g}(0) = 0.
\end{equation}
The stable Floquet bundle of $\Gamma$ is coupled with the BVP~\eqref{eq:po}--\eqref{eq:pobc} via the first variational equation. More precisely, we consider a second orbit segment $\mathbf{v_g}$, such that each point $\mathbf{v_g}(t)$ represents a vector associated with points $\mathbf{g}(t)$ of the orbit segment that solves~\eqref{eq:po}. The orbit segment $\mathbf{v_g}$ is a solution to the linearised flow such that $\mathbf{v_g}(0)$ is mapped to itself after one rotation around $\Gamma$. The length of $\mathbf{v_g}(0)$ is contracted after one rotation by the factor ${\rm exp}(T_\Gamma \, \lambda_{\rm s})$, which is the stable Floquet multiplier of $\Gamma$. We prefer formulating this in logarithmic form, which introduces the stable Floquet exponent $\lambda_{\rm s}$ to the first variational equation, rather than affecting the length of $\mathbf{v_g}(0)$. Therefore, the BVP~\eqref{eq:po}--\eqref{eq:pobc} is extended with the following system of equations:
\begin{equation}
\label{eq:FVE}
           \dot{\mathbf{v}}_{\mathbf{g}} = T_\Gamma \, \left[ {\rm D}_{\mathbf{g}} \mathbf{F}(\mathbf{g}) \, \mathbf{v_g} - \lambda_{\rm s} \, \mathbf{v_g} \right],
\end{equation}
\begin{equation}
\label{eq:FVEper}
  \mathbf{v_g}(1) - \mathbf{v_g}(0) = 0, 
\end{equation}
\begin{equation}
\label{eq:FVEbc}
   \mid\!\mid\! \mathbf{v_g}(0) \!\mid\!\mid = 1. 
\end{equation}
In particular, $\mathbf{v_g}(0) = \mathbf{v_g}(1)$ is the normalised vector that spans the local linearised isochron associated with $\mathbf{g}(0)$. 

We have not specified a phase condition and, indeed, we allow $\mathbf{g}$ to shift and start at any point $\gamma_\vartheta \in \Gamma$. Consequently, the linear bundle $\mathbf{v_g}$ will also shift such that $\mathbf{v_g}(0)$ still spans the local linearised isochron associated with $\mathbf{g}(0)$.

Phase shifting the periodic orbit and its linear bundle by continuation in this way has been performed before~\cite{ermentrout12}. However, the implementation in~\cite{ermentrout12} requires accurate knowledge of the coordinates of the point $\gamma_\vartheta$ in order to decide when to stop shifting. Our approach uses another BVP set-up to monitor the phase shift, so that both $\gamma_\vartheta$ and $\bv_\vartheta$ are determined up to {\sc Auto} accuracy. To this end, we introduce a third orbit segment $\bw$ that lies along $\Gamma$, with initial point $\bw(0)$ equal to $\mathbf{g}(0)$, and end point $\bw(1)$ equal to $\gamma_0$. The total integration time associated with this orbit segment $\bw$ is the fraction of the period $T_\Gamma$ that $\mathbf{g}(0)$ lies away from $\gamma_0$ along $\Gamma$; hence, it is directly related to the phase of $\mathbf{g}(0)$. We extend the BVP~\eqref{eq:po}--\eqref{eq:FVEbc} with the following system of equations: 
\begin{equation}
\label{eq:shift}
  \dot{\bw} = \nu \,  T_\Gamma \, \mathbf{F}(\bw),
\end{equation}
\begin{equation}
\label{eq:atphi}
  \bw(0) = \mathbf{g}(0), 
\end{equation}
\begin{equation}
\label{eq:atG0}
  \left[ \bw(1) - \gamma_0 \right] \cdot \bv_0^\perp = 0.
\end{equation}
Here, we do not impose $\bw(1) = \gamma_0$. Instead, condition~\eqref{eq:atG0} allows $\bw(1)$ to move in the linearisation of $\iso{0}$ at $\gamma_0$; this relaxation is necessary to ensure that the BVP remains well posed and the discretised problem has a solution. In practice, since $\bw(0) \in \Gamma$, the difference between $\bw(1)$ and $\gamma_0$ will be of the same order as the overall accuracy of the computation. Note that it is important to ensure $\nu \geq 0$ in equation~\eqref{eq:shift}, because $\bw(1)$ may diverge from $\gamma_0$ along $E(\gamma_0)$ otherwise. We found it convenient to start the calculation with $\nu = 1$, which corresponds to the orbit segment $\bw = \mathbf{g}$. 

The combined solution $\{ \mathbf{g}, \mathbf{v_g}, \bw \}$ to the multi-segment BVP~\eqref{eq:po}--\eqref{eq:atG0} represents a rotated version of $\Gamma$ and its stable Floquet bundle so that the head point is $\gamma_\varphi$ with phase $\varphi = 1 - \nu \, ({\rm mod} \ 1)$. We remark here that this extended set-up can also be used to compute $\iso{\varphi}$, for any phase $0 < \varphi < 1$, with the method  for $\iso{0}$ described in~\sref{sec:settings}; such a computation would approximate each isochron up to the same accuracy, without introducing an additional interpolation error.

%
\subsection{Continuation set-up for the phase reset}
\label{sec:2PBVP}
Recall the set-up for computing a phase reset by moving $\bu(0)$ transverse to $\iso{\vartheta}$, so that the end point $\bu(1)$ will move and lie on the linearisation of an isochron $\iso{\varphi}$ with a different phase $\varphi$. Here, the orbit segment $\bu$ is a solution of 
\begin{equation}
\label{eq:svf}
  \dot{\bu} = k \, T_\Gamma \, \mathbf{F}(\bu),
\end{equation}
for some $k \in \mathbb{N}$. The end point $\bu(1)$ should lie close to $\Gamma$ on the linearisation of $\iso{\varphi}$, for some $\varphi \in [0, 1)$. We stipulate that the rotated version of $\Gamma$ is shifted such that $\bu(1)$ lies close to $\mathbf{g}(0)$ along the direction $\mathbf{v_g}(0)$. Hence, we require the two boundary conditions
\begin{equation}
\label{eq:scon1}
  \left[ \bu(1) - \mathbf{g}(0) \right] \cdot \mathbf{v_g}(0)^{\phantom{\perp}} = \eta, 
\end{equation}
\begin{equation}
\label{eq:scon2}
  \left[ \bu(1) - \mathbf{g}(0) \right] \cdot \mathbf{v_g}(0)^{\perp} = 0, 
\end{equation}
where $\mathbf{v_g}(0)^{\perp}$ is the vector perpendicular to $\mathbf{v_g}(0)$. Here, $\eta$ measures the (signed) distance between $\bu(1)$ and $\mathbf{g}(0)$, which is along $\mathbf{v_g}(0)$. Since $\bu$ is a solution of~\eqref{eq:svf} and $k \in \mathbb{N}$, the initial point $\bu(0)$ has the same phase as the last point $\bu(1)$, and the combined multi-segment BVP~\eqref{eq:po}--\eqref{eq:scon2} ensures that $\bu(1)$ has (approximate) phase $1 - \nu \, ({\rm mod} \ 1)$. In practice, we should choose $k \in \mathbb{N}$ large enough such that $\eta < \eta_{\rm max}$. If $\bu(0)$ lies close to $\Gamma$, it will be sufficient to set $k = 1$. In order to consider phase resets of large perturbations, for which $\bu(0)$ starts relatively far away, we need $k > 1$, to allow for sufficient time to let $\bu$ converge and have $\bu(1)$ lie close to $\Gamma$.

At this stage, the multi-segment BVP~\eqref{eq:po}--\eqref{eq:scon2} is a system of ${\sf NDIM} = 8$ ordinary differential equations (for the case of a planar system), with ${\sf NBC} = 10$ boundary conditions, and ${\sf NPAR} = 4$ free parameters, namely, $T_\Gamma$, $\lambda_{\rm s}$, $\nu$, and $\eta$; the period $T_\Gamma$ and stable Floquet exponent $\lambda_{\rm s}$ must remain free parameters to ensure that the discretised problem has a solution, but their variation will be almost zero. Hence, ${\sf NDIM} - {\sf NBC} + {\sf NPAR} = 2 \neq 1$, and one more condition is needed to obtain a one-parameter family of solutions.

The final step in the set-up is to impose an extra condition that specifies how $\bu(0)$ moves along an arc or closed curve in the phase plane. Consequently, since $k \, T_\Gamma$ is fixed, the orbit segment $\bu$ changes, so that $\bu(1)$ will move as well, and $\mathbf{g}(0)$, along with $\mathbf{v_g}(0)$ will shift accordingly. This causes a variation in $\nu$ to maintain $\bw(0) = \mathbf{g}(0)$, and these $\nu$-values precisely define the phase-response curve in the continuation run, where the position along the chosen arc or closed curve is the argument.

To compute the PRC, we need to let $\bu(0)$ traverse the closed curve $\{\gamma_\vartheta + A \, \mathbf{d} \;|\; \vartheta \in [0, 1) \}$ obtained by the (instantaneous) perturbation of $\Gamma$ in the direction $\mathbf{d}$ for distance $A$. We can impose this relatively complicated path on $\bu(0)$ by including another system of equations to the multi-segment BVP, namely, the BVP that defines $\Gamma$ in terms of another rotated orbit segment $\mathbf{g_u}$. Furthermore, in order to keep track of the phase $\vartheta$ along this path, we introduce another segment $\bw_{\bu}$ that plays the same role as $\bw$ in \sref{sec:rotate}; compare with equations~\eqref{eq:po}--\eqref{eq:pobc} and~\eqref{eq:shift}--\eqref{eq:atG0}. In other words, we extend the BVP~\eqref{eq:po}--\eqref{eq:scon2} by the following system of equations
\begin{equation}
\label{eq:gu}
  \dot{\mathbf{g}}_{\mathbf{u}} = \widehat{T}_\Gamma \, \mathbf{F}(\mathbf{g_u}),
\end{equation}
\begin{equation}
\label{eq:gubc}
   \mathbf{g_u}(1) - \mathbf{g_u}(0) = 0.
\end{equation}
\begin{equation}
\label{eq:wu}
  \dot{\bw}_{\mathbf{u}} = (1 - \vartheta) \, \widehat{T}_\Gamma \, \mathbf{F}(\bw_{\bu}),
\end{equation}
\begin{equation}
\label{eq:wuatphi}
  \bw_{\bu}(0) = \mathbf{g_u}(0), 
\end{equation}
\begin{equation}
\label{eq:wuatG0}
  \left[ \bw_{\bu}(1) - \gamma_0 \right] \cdot \bv_0^\perp = 0.
\end{equation}
Here, we decrease $\vartheta$ from $1$ to $0$, during which $\bw_{\bu}$ grows and $\mathbf{g_u}$ tracks $\gamma_\vartheta$. In order for a solution to exist, the periods $T_\Gamma$ and $\widehat{T}_\Gamma$ must be two different free parameters, although they remain constant (and equal) to within the accuracy of the computation. The phase reset is now obtained by imposing
\begin{equation}
\label{eq:reset}
  \bu(0) = \mathbf{g_u}(0) + A \, \mathbf{d}.
\end{equation}
The multi-segment BVP~\eqref{eq:po}--\eqref{eq:reset} is now a system of dimension ${\sf NDIM} = 12$, with
${\sf NBC} = 17$ boundary conditions, and ${\sf NPAR} = 6$ free parameters, which are $T_\Gamma$, $\lambda_{\rm s}$, $\nu$,  $\eta$, $\widehat{T}_\Gamma$, and either $\vartheta$ or $A$. Since, ${\sf NDIM} - {\sf NBC} + {\sf NPAR} = 1$, we obtain a one-parameter solution family by continuation. As the first solution in the continuation, we use the known solution $\mathbf{g} = \bw = \mathbf{u} =\mathbf{g_u} = \Gamma$, which starts with the head point $\gamma_0$, the associated stable linear bundle $\bv_0$ that we assumed has been pre-computed, and $\bw_{\bu} = \gamma_0$; then $T_\Gamma = \widehat{T}_\Gamma$ and $\lambda_{\rm s}$ are set to their known computed values, $\eta = A = 0$, and $\nu = \vartheta = 1$. Initially, $k = 1$, and one should monitor $\eta$ to make sure it does not exceed $\eta_{\rm max}$.

To obtain the PTC or PRC we first perform a homotopy step, where we fix $\vartheta = 1$ and vary the amplitude $A$ until the required value is reached. This continuation run produces a one-parameter family of solutions representing the effect of a reset of varying amplitude $A$ from the point $\gamma_0$. In the main continuation run, we then fix $A$ and decrease $\vartheta$ until $\vartheta = 0$, so that it covers the unit interval; the associated solution family of the multi-segment BVP~\eqref{eq:po}--\eqref{eq:reset}, hence, provides the resulting phase shift $\vartheta_{\rm new} := 1 - \nu \, ({\rm mod} \ 1)$ as a function of the phase $\vartheta_{\rm old} := \vartheta$ along the perturbed periodic orbit.

%
\section{Illustration of the method with a model example}
\label{sec:winfree}
We illustrate our method for computing a PTC with a constructed example, namely, a parametrised version of the model introduced by Winfree~\cite[Chapter~6]{Winfree2001}, which we also used in~\cite{lko-siads2015}; it is given in polar coordinates as
\begin{displaymath}
  \left\{ \begin{array}{rcc}
      \dot{r} &=& (1-r) \, (r-a) \, r,  \\
      \dot{\psi} &=& -1 - \omega \, (1-r).
  \end{array} \right.
\end{displaymath}
In Euclidean coordinates, the system becomes
\begin{equation}
\label{eq:WinEx}
  \left\{ \begin{array}{rcccl}
      \dot{x} &=&  (1 - \sqrt{x^2 + y^2}) \, \left(  x \, (\sqrt{x^2 + y^2} - a) + \omega y \right) + y, \\[2mm]
      \dot{y} &=&  (1 - \sqrt{x^2 + y^2}) \, \left(  y \, (\sqrt{x^2 + y^2} - a)  - \omega x \right) - x.
  \end{array} \right.
\end{equation}
Note that this system is invariant under any rotation about the origin; moreover, its frequency of rotation only depends on $r = \sqrt{x^2 + y^2}$; see~\cite{lko-siads2015} for details. We now fix the parameters to $a = 0$ and $\omega = -0.5$, as in~\cite{lko-siads2015}. Then the unit circle is an attracting periodic orbit $\Gamma$ with period $T_\Gamma = 2\pi$ and the origin is an unstable equilibrium $\xeq$.

%
\subsection{Computing the PTC}
\label{sec:ptcexample}
We choose $\gamma_0 = (1, 0)$ and compute the normalised linear direction associated with its isochron as $\bv_0 \approx (-0.83, -0.55)$. As was explained in \sref{sec:2PBVP}, the computation is performed in two separate continuation runs: first, we apply a perturbation to the point $\gamma_0$ in a fixed direction $\mathbf{d}$, where we vary the amplitude $A$ from $0$ to $0.75$ during the homotopy step. Next, we fix $A = 0.75$ and apply the same perturbation to each point $\gamma_\vartheta \in \Gamma$. For the purpose of visualising the computational set-up, we choose the (somewhat unusual) direction $\mathbf{d} = (-1,0)$ and set the maximum distance along the linearised isochron to the relatively large value of $\eta_{\rm max} = 0.2$.

%
\begin{figure}[t!]
    \centering
    \includegraphics{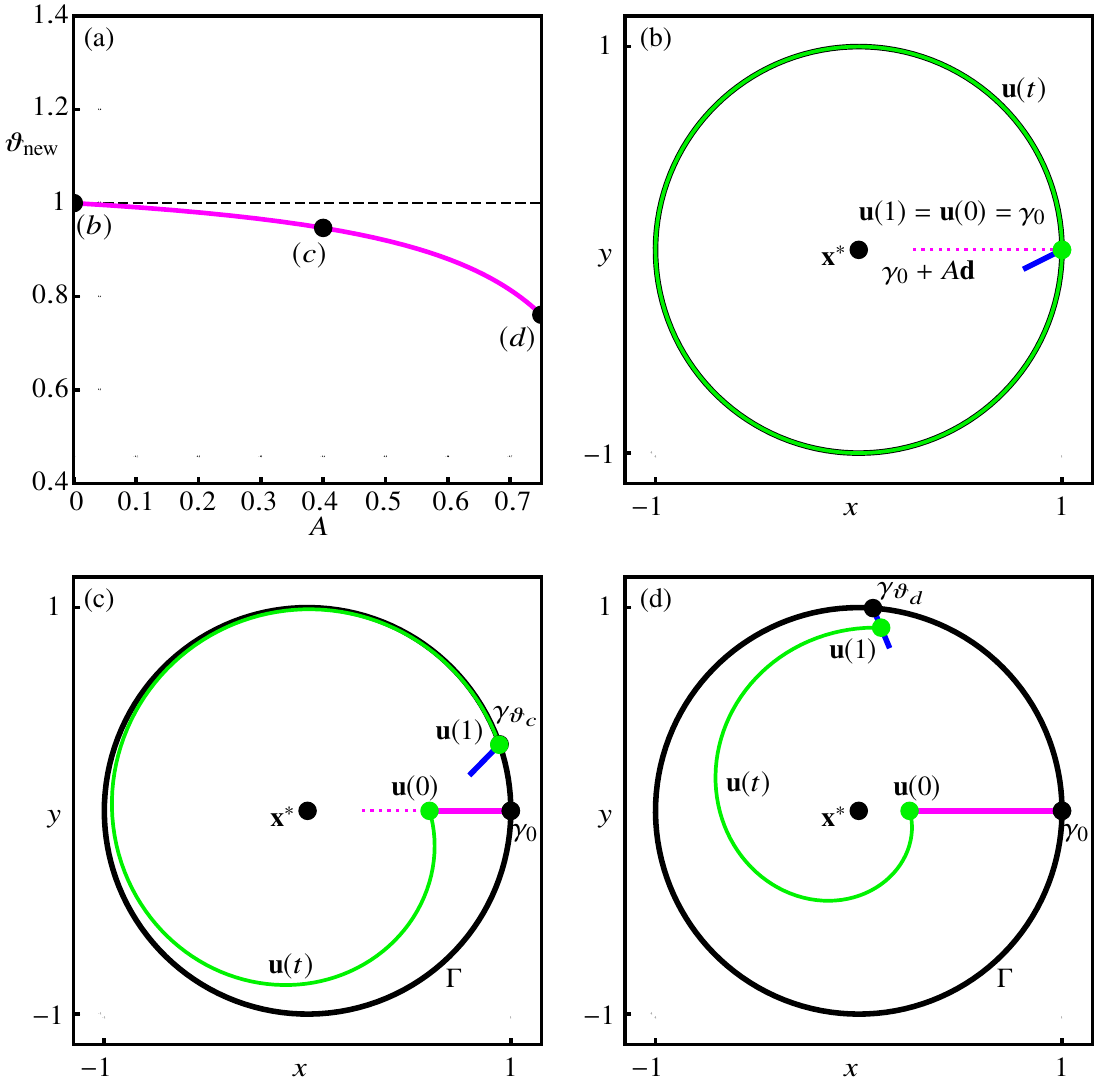}
    \caption{\label{fig:comp1}
    Phase reset of system~\eqref{eq:WinEx} at fixed $\gamma_0$ in the direction $\mathbf{d} = (-1, 0)$ with amplitude $A \in [0, 0.75]$~(a) and continuation set-up at the three labelled points~(b)--(d); here, panel~(b) shows the initial set-up when $A = 0$ and $\vartheta_{\rm new} = 1$, in panel~(c) the continuation has progressed to $A = 0.4$ and $\vartheta_{\rm new} = \vartheta_c \approx 0. 96$, and in panel~(d) $A = 0.75$ has been reached and $\vartheta_{\rm new} = \vartheta_d \approx 0.76$.}
\end{figure}
%
%
The first continuation run of the multi-segment BVP~\eqref{eq:po}--\eqref{eq:reset} is illustrated in \fref{fig:comp1}. Here, the free amplitude $A$ increases while $\vartheta = 1 = 0 \, ({\rm mod} \ 1)$ is fixed and, hence, the perturbation is always applied at $\gamma_0$ and grows in size. \Fref{fig:comp1}(a) shows the resulting phase $\vartheta_{\rm new}$ as a function of $A$. Three points are labelled, indicating the three stages during the continuation that are illustrated in panels~(b), (c) and~(d). In each of these panels we show the periodic orbit $\Gamma$ in black, and the current orbit segment $\bu$ of the continuation run in green. Note that $\Gamma$ is rotated here and its head point $\mathbf{g}(0)$ lies at the point on $\Gamma$ with phase $\vartheta_{\rm new}$. A short segment of the associated linearisation of the isochron of $\gamma_{\vartheta_{\rm new}}$ is shown in blue. We do not plot the orbit segment $\bw$ that determines the value of $\vartheta_{\rm new}$, but it follows $\Gamma$ from $\mathbf{g}(0)$ back to $\mathbf{g}(0)$ and then extends (approximately) along $\Gamma$ to $\gamma_0$. Indeed, notice in \fref{fig:comp1}(a) that $\vartheta_{\rm new}$ is decreasing, which means that $\nu > 1$ is increasing so that $\bw$ becomes longer. We also do not show the orbit segments $\mathbf{g_u}$ and $\bw_{\bu}$ that determine the phase $\vartheta = \vartheta_{\rm old}$ at which the perturbation is applied, because $\vartheta_{\rm old} = 1$ is fixed in this continuation run.

\Fref{fig:comp1}(b) shows the initial set-up, with $\mathbf{g} = \bw = \mathbf{u} = \mathbf{g_u} = \Gamma$, $\bw_{\bu} = \gamma_0$, $T_\Gamma = \widehat{T}_\Gamma$ and $\lambda_{\rm s}$ set to their known values, and $\nu = 1$, $\eta = A = 0$, with $k = 1$ and $\vartheta = 1$. The dotted line segment in \fref{fig:comp1}(b) indicates the direction $\mathbf{d}$ of the intended perturbation away from $\gamma_0$; its length is the maximal intended amplitude $A = 0.75$. An intermediate continuation step when $A = 0.4$ is shown in \fref{fig:comp1}(c). The perturbation has pushed $\bu(0)$ out along $\mathbf{d}$, such that $\bu(1)$ now lies (approximately) on the linearised isochron, parametrised as $\mathbf{g}(0) + \eta \, \mathbf{v_g}(0)$ with $0 < \eta \leq \eta_{\rm max}$, associated with the rotated head point $\mathbf{g}(0) = \gamma_{\vartheta_c}$, where $\vartheta_c \approx 0.96$. Note that the orbit segment $\bw$ (not shown) has now changed from its initialisation to match the solution to subsystem~\eqref{eq:shift}--\eqref{eq:atG0} with $\nu \approx 1.04$. \Fref{fig:comp1}(d) illustrates the last step of the first continuation run, when $A = 0.75$. The head point $\mathbf{g}(0) \in \Gamma$ has rotated further to $\gamma_{\vartheta_d}$ with $\vartheta_d = 1 - \nu \approx -0.24 = 0.76 \, ({\rm mod} \ 1)$. Notice that $\bu(1)$ lies quite far along the linearised isochron, because we allow a relatively large distance $\eta$. The corresponding orbit segment $\bu$ is determined for an integration time of only one period, that is, for $k = 1$. We show this case for illustration purposes, but in practice, it would be worth choosing a smaller value for $\eta_{\rm max}$, so that $\bu$ would be extended, and the integer multiple of $T_\Gamma$ set to $k = 2$, before reaching $A = 0.75$.

The second continuation run uses the fixed perturbation of size $A = 0.75$ along $\mathbf{d} = (-1, 0)$, and varies the phase $\vartheta$ at which it is applied. Since $\vartheta$ controls the integration time associated with the orbit segment $\bw_{\bu}$, the multi-segment BVP~\eqref{eq:gu}--\eqref{eq:wuatG0} with solution $\{ \mathbf{g_u}, \bw_{\bu} \}$ and parameter $\widehat{T}_\Gamma$ plays an important role now. For each $\vartheta$, the head point $\mathbf{g_u}(0)$ of $\mathbf{g_u}$ lies (approximately) at $\gamma_\vartheta \in \Gamma$, and $\bw_{\bu}$ represents the remaining part of $\Gamma$ from $\gamma_\vartheta$ to $\gamma_0$; hence, the total integration time of $\bw_{\bu}$ is the fraction $1 - \vartheta$ of $\widehat{T}_\Gamma$, which is equal, up to the computational accuracy, to the period $T_\Gamma$ of $\Gamma$.

%
\begin{figure}[t!]
    \centering
    \includegraphics{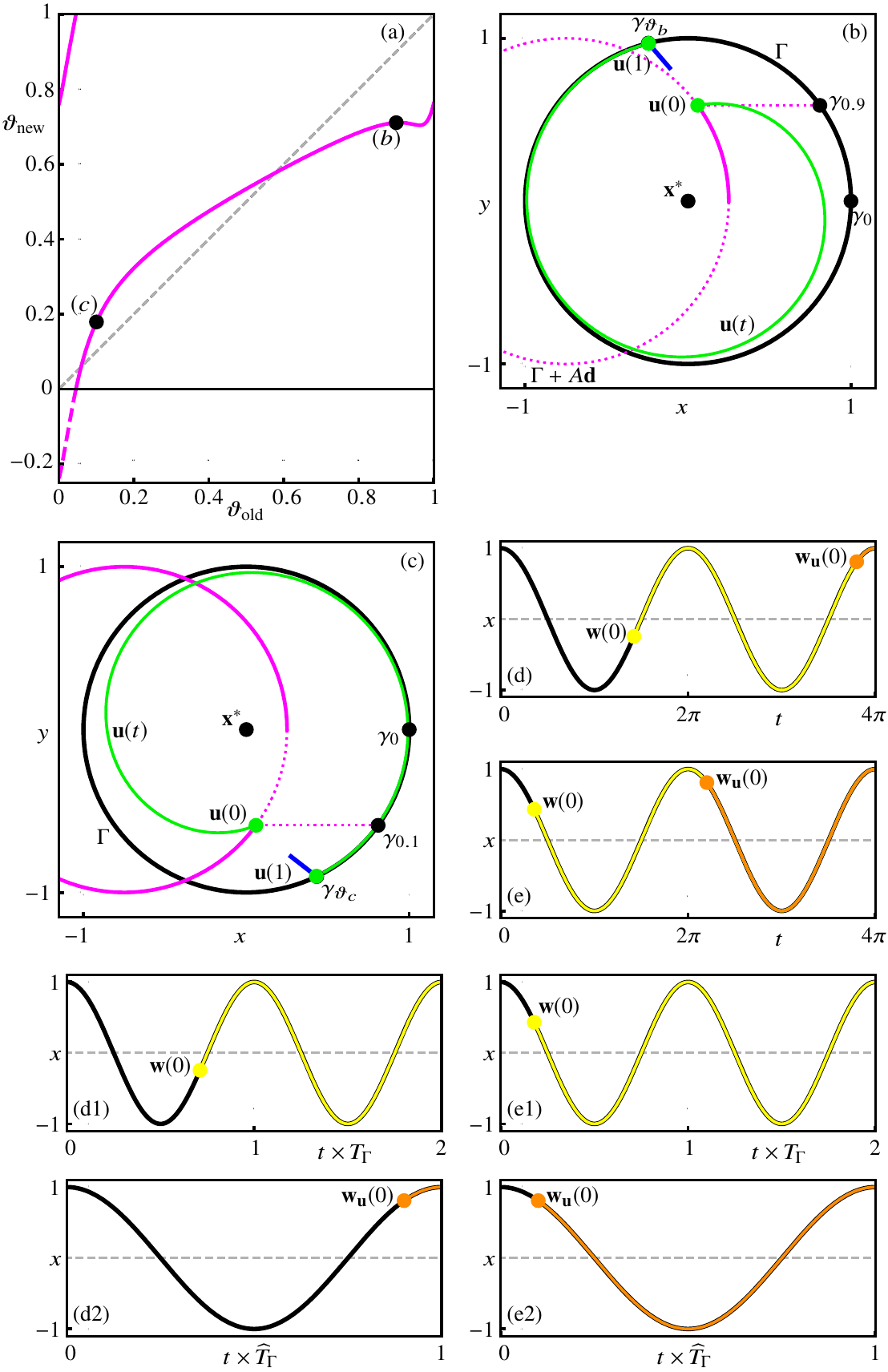}
    \caption{\label{fig:comp2}
    PTC of $\Gamma$ in system~\eqref{eq:WinEx} for $\mathbf{d} = (-1, 0)$ and $A = 0.75$~(a), and continuation set-up at $\vartheta_{\rm old} = 0.9$~(b) and at $\vartheta_{\rm old} = 0.1$~(c) with $\bw$ and $\bw_{\bu}$ in~(d), (d1), (d2) and~(e), (e1), (e2), respectively.}
\end{figure}
\Fref{fig:comp2} illustrates different aspects of this continuation run. As $\vartheta_{\rm old} = \vartheta$ decreases from $1$, the multi-segment BVP~\eqref{eq:po}--\eqref{eq:reset} determines the orbit segment $\bu$ with $\bu(0) = \gamma_\vartheta + A \, \mathbf{d}$ and uses the rotated orbit segment $\mathbf{g}$ and $\bw$ to establish the resulting phase $\vartheta_{\rm new} = 1 - \nu \, ({\rm mod} \ 1)$ of $\bu(1)$. Panel~(a) shows the PTC computed for $A = 0.75$. Note that $\nu$ takes values in the covering space $\mR$; the output is then folded onto the unit torus by taking $\vartheta_{\rm new} = 1 - \nu \, ({\rm mod} \ 1)$, giving the solid curve in \fref{fig:comp2}. The points labelled~(b) and~(c) in this panel correspond to $\vartheta_{\rm old} = 0.9$ and $\vartheta_{\rm old} = 0.1$, respectively. The continuation set-up for these two cases is shown in the corresponding panels~(b) and~(c). As in \fref{fig:comp1}, the periodic orbit $\Gamma$ is black and $\bu$ is green. The path traced by the initial point $\bu(0)$ is the magenta dotted circle, which is $\Gamma$ shifted by $A = 0.75$ in the direction $\mathbf{d} = (-1, 0)$; hence, for fixed $\vartheta$, the point $\bu(0)$ corresponds to the perturbation of the point $\gamma_\vartheta \in \Gamma$ that lies horizontally to the right of $\bu(0)$, as indicated by the magenta dotted line segment. The end point $\bu(1)$ lies on the linearised isochron, parametrised as $\mathbf{g}(0) + \eta \, \mathbf{v_g}(0)$ with $0 < \eta \leq \eta_{\rm max}$, associated with the rotated head point of $\mathbf{g}$, which is determined by subsystem~\eqref{eq:po}--\eqref{eq:FVEbc}. The phase of this head point is given by $\vartheta_{\rm new} = 1 - \nu \, ({\rm mod} \ 1)$, where $\nu$ is determined from subsystem~\eqref{eq:shift}--\eqref{eq:atG0} that defines the orbit segment $\bw$. 

Hence, the two orbit segments $\bw$ and $\bw_{\bu}$ essentially determine the PTC, that is, the map $P: \vartheta_{\rm old} \mapsto \vartheta_{\rm new}$. Their $x$-coordinate is plotted versus time in panel~(d) for $\vartheta_{\rm old} = 0.9$ and in panel~(e) for $\vartheta_{\rm old} = 0.1$, respectively, overlaid on two copies of $\Gamma$ (black curve), that is, time $t$ runs from $0$ to $4 \pi$. The further panels~(d1) and~(d2) for $\vartheta_{\rm old} = 0.9$ and in panels~(e1) and~(e2) for $\vartheta_{\rm old} = 0.1$ show $\bw$ (yellow curve)  and $\bw_{\bu}$ (orange curve) individually over the fraction of the periods of $T_\Gamma$ and $\widehat{T}_\Gamma$, respectively. Note that both $\bw$ and $\bw_{\bu}$ end at $t = 4 \pi$ and $x = 1$, as required, but their initial points differ. As $\vartheta$ decreases from $1$ to $0$ during the continuation, the orbit segment $\bw_{\bu}$ lengthens as expected, but note that $\bw$ lengthens as well; this is due to the (near-)monotonically increasing nature of the PTC for this example.

%
\subsection{Loss of invertibility}
\label{sec:noinvert}

Recall that any PTC for $A =0$ is the identity, and it is invertible for sufficiently small amplitude $A$ of the perturbation, because its graph remains a 1:1 torus knot on the torus parametrised by the two periodic variables $\vartheta_{\rm old}$ and $\vartheta_{\rm new}$. However, the PTC in \fref{fig:comp2}(a) for $A = 0.75$ is no longer near the identity: it is not injective and, hence, not invertible. 

To show how injectivity of the PTC is lost as $A$ is increased, we consider again model~\eqref{eq:WinEx}, but now with $a = 0.25$; see also~\cite{lko-siads2015}. Apart from the attracting unit circle $\Gamma_{\rm s} = \Gamma$ with period $T_\Gamma = 2 \pi$, there exists then also a repelling circle $\Gamma_{\rm u}$ with radius $r = a = 0.25$ and period $2 \pi / (1 + \omega \, (1 - a)) = 3.2 \, \pi$; note that $\Gamma_{\rm u}$ forms the boundary of the basins of attraction of both $\Gamma_{\rm s}$ and the equilibrium $\xeq$ at the origin, which is now attracting. 

%
\begin{figure}[t!]
    \centering
    \includegraphics{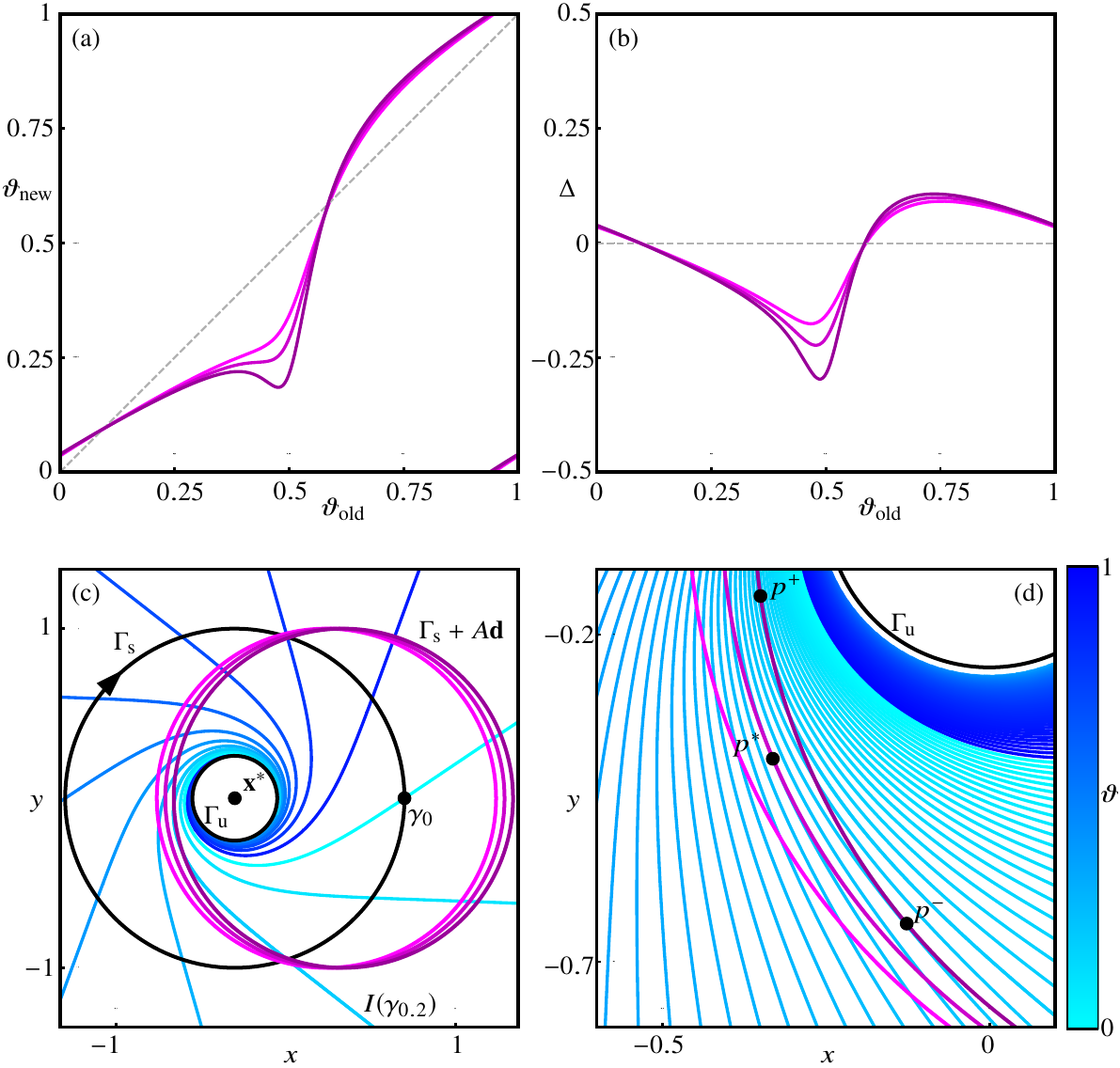}
    \caption{\label{fig:simple}
      Phase resets of $\Gamma_{\rm s}$ in system~\eqref{eq:WinEx} with $a=0.25$ for amplitudes $A \in \{0.54, 0.59, 0.64\}$ (increasingly darker shades of magenta). The three PTCs are shown in panel~(a) and the corresponding PRCs in panel~(b). The three perturbed cycles are shown in panel~(c) together with $\Gamma_{\rm s}$ and ten of its isochrons that are uniformly distributed over one period; the enlargement near $\Gamma_{\rm u}$ in panel~(d) shows them with 100 uniformly distributed isochrons of $\Gamma_{\rm s}$ and points of tangency at $p^*$ and $p^\pm$. Isochrons are coloured according to the colour bar.}
\end{figure}
%
%
We consider three resets of $\Gamma_{\rm s}$ of the form $\Gamma_{\rm s} + A \, \mathbf{d}$ in the positive direction $\mathbf{d} = (1, 0)$ and with $A = 0.54$, $A = 0.59$, and $A = 0.64$. \Fref{fig:simple} shows the three corresponding PTCs, the corresponding PRCs, and the perturbed cycles $\Gamma_{\rm s} + A \, \mathbf{d}$ in increasingly darker shades of magenta as $A$ increases in panels~(a), (b), and~(c), respectively. Panel~(a) shows that the first PTC for $A = 0.54$ is injective and invertible. As $A$ is increased to approximately $A = 0.59$, the graph has a cubic tangency near $(\vartheta_{\rm old}, \vartheta_{\rm new}) = (0.45, 0.24)$, because the associated map $P$ has an inflection point at $\vartheta_{\rm old} \approx 0.45$. For larger values of $A$, such as for $A = 0.64$, the PTC has a local maximum followed by a local minimum and is, hence, no longer invertible. Note from \fref{fig:simple}(b) that this qualitative change of the PTC does not lead to a corresponding qualitative change of the PRC. 

\Fref{fig:simple}(c) and the enlargement near the basin boundary $\Gamma_{\rm u}$ in panel~(d) show that the loss of injectivity of the PTC is due to a cubic tangency between the perturbed cycle $\Gamma_{\rm s} + A \mathbf{d}$ and the foliation of the basin of $\Gamma_{\rm s}$ by (forward-time) isochrons; ten isochrons are shown in panel~(c) and one hundred in panel~(d), distributed uniformly in phase and coloured according to the colour bar. The left-most light-magenta cycle for $A = 0.54$ is transverse to all isochrons. The middle magenta cycle for $A=0.59$, on the other hand, has a single cubic tangency (approximately) with the isochron $\iso{0.24}$ of phase $\vartheta = 0.24$ at the point $p^* \approx (-0.33, -0.39)$, shown in panel~(d). For larger $A$, as for the right-most dark-magenta cycle for $A=0.64$, there are now two quadratic tangencies with two different isochrons, namely, (approximately) with $\iso{0.22}$ and $\iso{0.19}$ at the points $p^+ \approx (-0.35,-0.14)$ and $p^- \approx (-0.13,-0.64)$, respectively. As a result, all isochrons that intersect the perturbed cycle between $p^+$ and $p^-$ intersect three times; hence, the map $P$ from $\vartheta_{\rm old}$ to $\vartheta_{\rm new}$ is no longer invertible. Note that $p^+$ and $p^-$ correspond to the local maximum and local minimum of the PTC in panel~(a), respectively.

%
\section{Phase resetting in the FitzHugh-Nagumo model}
\label{sec:fhn}
We now illustrate the capability of our method by computing the phase response of a periodic orbit in the FitzHugh--Nagumo system~\cite{fitzhugh61, Nagumo1962}. This model is an iconic example that motivated early work on isochrons and phase response curves; in particular, it has a very complicated geometry of isochrons with regions of extreme phase sensitivity ~\cite{lko-chaos2014, Winfree2001}. The FitzHugh--Nagumo system is given by the equations
\begin{equation}
\label{eq:FHN}
  \left\{ \begin{array}{rcl}
      \dot{x} &=& {\displaystyle c \, \left( y + x - \tfrac{1}{3} \, x^3 + z \right)}, \\[2mm]
      \dot{y} &=& {\displaystyle -\frac{x - a + by}{c}}.
  \end{array} \right.
\end{equation}
We set $a = 0.7$, $b = 0.8$, and $z = -0.4$, as in~\cite{Winfree2001}, and fix $c = 2.5$, as was done in~\cite{lko-siads2015}. For these parameter values, there exists an attracting periodic orbit $\Gamma$ with period $T_\Gamma \approx 10.71$ and a repelling equilibrium $\xeq \approx (0.9066, -0.2582)$. The parameter $c$ is a time-scale parameter, the increase of which makes the $x$-variable faster than the $y$-variable. It plays an important role in the onset of phase sensitivity due to an accumulation of isochrons in a narrow region close to the slow manifold~\cite{lko-chaos2014}, which is associated with the occurrence of sharp turns in the isochrons of $\Gamma$; see also~\cite{om-siads2010}. For the chosen value of $c = 2.5$, one finds both strong phase sensitivity and sharp turns, which makes the computation of any phase response quite challenging.

%
\begin{figure}[t!]
  \centering
  \includegraphics{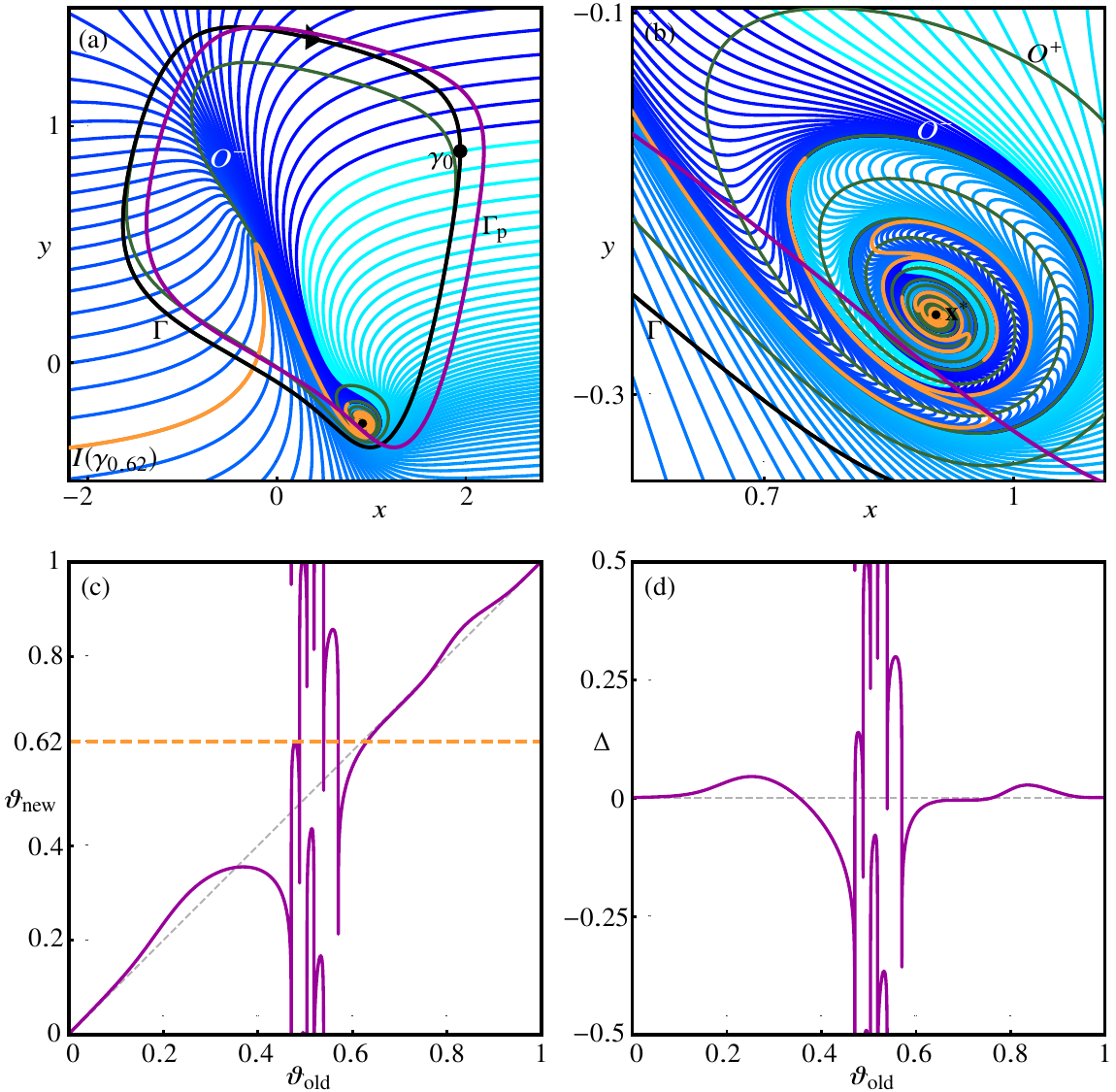}
  \caption{\label{fig:PRC}
    Phase reset for the FitzHugh--Nagumo model~\eqref{eq:FHN}. Panel~(a) shows the periodic orbit $\Gamma$ (black), the perturbed cycle $\Gamma + A \mathbf{d}$  (magenta) with $\mathbf{d} = (1,0)$ and $A = 0.25$, two trajectories $O^+$ and $O^-$ (green), and 100 isochrons uniformly distributed in phase; isochrons are coloured according to the colour bar in \fref{fig:simple}, and the isochron $I(\gamma_{0.62})$ is highlighted in orange. Panel~(b) is an enlargement near the equilibrium $\xeq$, and panels~(c) and~(d) show the corresponding PTC and PRC, respectively; the dashed orange line in panel~(c) indicates the phase of $I(\gamma_{0.62})$.}
\end{figure} 
%
%
\Fref{fig:PRC} illustrates the phase reset for the FitzHugh--Nagumo model~\eqref{eq:FHN} after a perturbation in the $x$-direction $\mathbf{d} = (1,0)$ of amplitude $A = 0.25$. Panel~(a) and the enlargement near the equilibrium $\xeq$ in panel~(b) show how the perturbed cycle $\Gamma + A \, \mathbf{d}$ intersects the isochrons of $\Gamma$, of which 100 are shown equally distributed in phase and coloured according to the colour bar in \fref{fig:simple}. In particular, one notices quite a few instances in panel~(b) of quadratic tangencies between the perturbed cycle and different isochrons; one such isochron is the highlighted $I(\gamma_{0.62})$. The green curves $O^+$ and $O^-$ in panels~(a) and~(b) are two special trajectories, along which the foliation by forward-time isochrons of $\Gamma$ has quadratic tangencies with the foliation by backward-time isochrons (not shown) of the focus $\xeq$. Tangencies between these two foliations were introduced in~\cite{lko-siads2015}, where we argued that such tangencies give rise to sharp turns of isochrons. We remark that the two trajectories $O^+$ and $O^-$ of quadratic tangencies arises at a specific value $c^* < 2.5$ where one finds a cubic tangency between the two foliations, called a cubic isochron foliation tangency or CIFT for short; see~\cite{lko-siads2015} for details. The relevance of the special trajectories $O^+$ and $O^-$ in the present context is that along them the isochrons of $\Gamma$ have sharp turns as they approach $\xeq$. This can clearly be seen in \fref{fig:PRC}(b); as the highlighted isochron $I(\gamma_{0.62})$ illustrates, the turns along $O^-$ are so sharp that $I(\gamma_{0.62})$ appears to retrace itself along certain segments. Since this happens for all isochrons of $\Gamma$, one finds extreme phase sensitivity near the trajectory $O^-$. Moreover, quadratic tangencies of the perturbed cycle with isochrons of $\Gamma$ occur near both $O^+$ and $O^-$. Hence, the number of intersection of $\Gamma + A \mathbf{d}$ with $O^+$ and $O^-$ gives an indication of how many quadratic tangencies the perturbed cycle has with different isochrons.

As we have seen in \sref{sec:noinvert}, any such quadratic tangency between $\Gamma + A \, \mathbf{d}$ and an isochron is associated with a local maximum or minimum of the PTC, which is, therefore, not expected to be invertible. \Fref{fig:PRC}(c) presents the PTC computed with our method as a continuous curve shown on the $(\vartheta_{\rm old}, \vartheta_{\rm new})$ unit torus. Clearly, its graph is quite intriguing and features six local maxima and six local minima. Observe that the local maxima correspond to quadratic tangencies near $O^+$, while the sharper local minima correspond to quadratic tangencies near $O^-$; in particular, the tangency with the highlighted isochron $I(\gamma_{0.62})$ near $O^+$ in panel~(b) gives rise to a local maximum of the PTC in panel~(c), where the graph has a tangency with the dashed orange line at $\vartheta_{\rm new} = 0.62$. Notice that $I(\gamma_{0.62})$ intersects the perturbed cycle $\Gamma + A \mathbf{d}$ in panel~(b), and hence, the PTC in panel~(c), five more times. The associated PRC of the change in phase $\Delta = \vartheta_{\rm new} - \vartheta_{\rm old}$ is shown in panel~(d); it is also quite a complicated curve with corresponding local maxima and minima. The PTC and PRC both have six near-vertical segments at $\vartheta_{\rm old} \approx 0.47$, $0.49$, $0.50$, $0.52$, $0.54$, and $0.57$; such large gradients arise near the local minima because of the extreme phase sensitivity near $O^-$. 

%
\begin{figure}[t!]
  \centering 
  \includegraphics{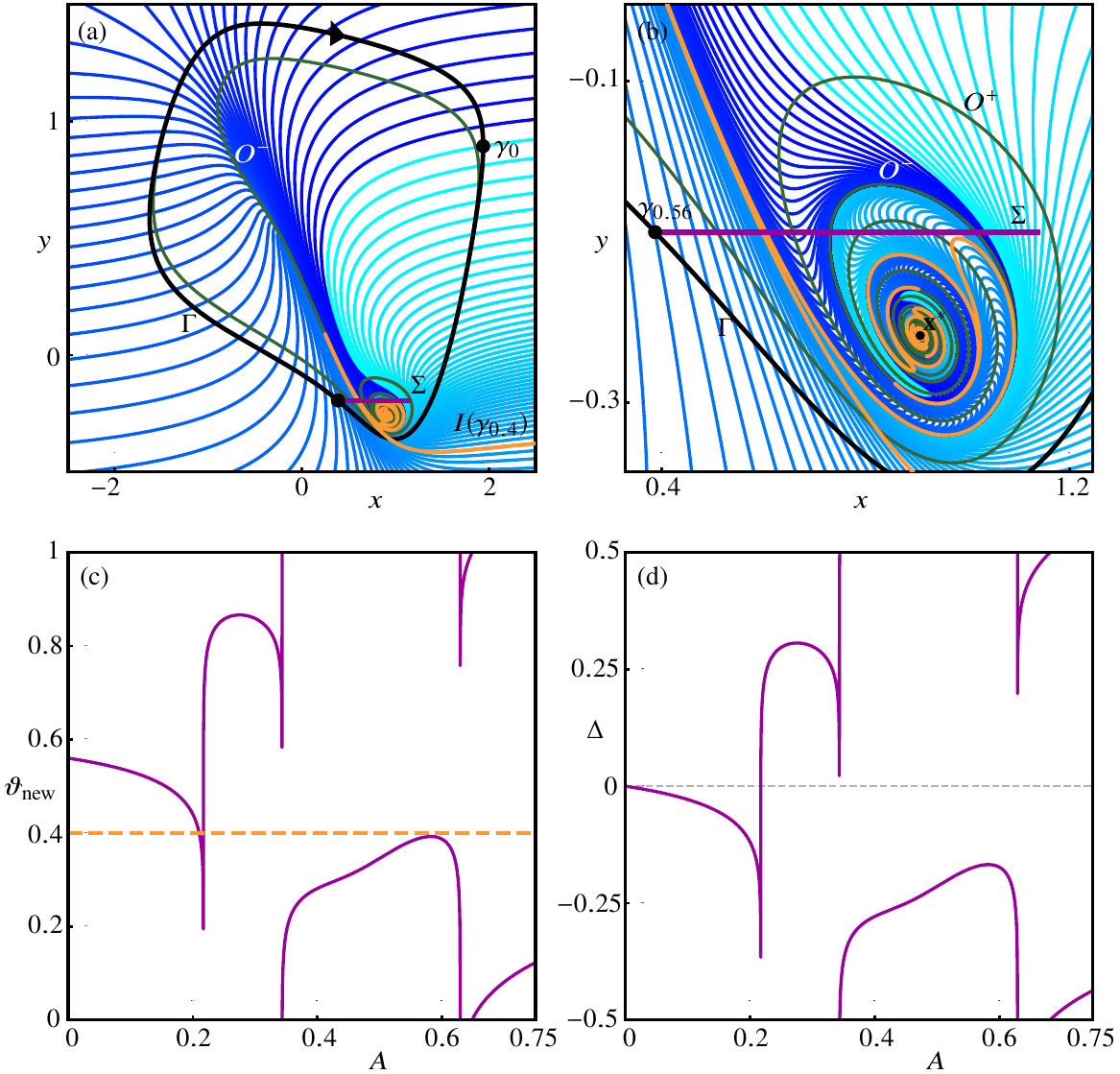}
  \caption{\label{fig:ARC} 
    Phase response along the line segment $\gamma_{0.56} + A \, \mathbf{d}$ with $\mathbf{d} = (1,0)$ and $A \in \left[ 0 , 0.75 \right]$ in the FitzHugh--Nagumo model~\eqref{eq:FHN}. Panel~(a) shows the periodic orbit $\Gamma$ (black), the line segment of perturbations (magenta) starting at point $\gamma_{0.56} \in \Gamma$, the two trajectories $O^+$ and $O^-$ (green), and 100 isochrons uniformly distributed in phase; isochrons are coloured according to the colour bar in \fref{fig:simple}, and the isochron $I(\gamma_{0.63})$ is highlighted in orange. Panel~(b) is an enlargement near the equilibrium $\xeq$, and panels~(c) and~(d) show the periodic variables $\vartheta_{\rm new}$ and $\Delta = \vartheta_{\rm new} - \vartheta_{\rm old}$, respectively, as a function of $A$.}
\end{figure} 
%
%
\Fref{fig:ARC} illustrates our continuation approach for another type of resetting experiment, where phase and direction of the perturbation are fixed but its magnitude varies. Specifically, we calculate the asymptotic phase of points that are perturbed from $\gamma_{0.56} \in \Gamma$ in the positive $x$-direction $\mathbf{d} = (1,0)$ with amplitude $A \in [0 , 0.75]$. As panels~(a) and~(b) show, the corresponding line segment $\gamma_{0.56} + A \, \mathbf{d}$ passes through the phase-sensitive region of accumulating isochrons near $\xeq$, where it intersects $O^+$ and $O^-$ several times. To compute the phase response, we first rotate $\Gamma$ and, consequently, the entire multi-segment BVP~\eqref{eq:po}--\eqref{eq:reset}, such that the head point $\mathbf{g}(0)$ of $\Gamma$ lies at $\gamma_{0.56}$. We then proceed as in the first continuation run in \sref{sec:ptcexample} to obtain $\vartheta_{\rm new}$ as a function of $A$. The resulting phase responses of  $\vartheta_{\rm new}$ and $\Delta \vartheta_{\rm new}$ are shown in panels~(c) and~(d), respectively; note that $\Delta \vartheta_{\rm new}$ is obtained from $\vartheta_{\rm new}$ by a fixed shift of $\vartheta_{\rm old} = 0.56$. The resulting phase response as a function of $A$ also shows vertical segments near three local minima, which are again directly associated with the three points where the line segment $\gamma_{0.56} + A \, \mathbf{d}$ intersects the trajectory $O^-$. Notice that the turns of the isochrons along $O^-$ are so very sharp that one will find a quadratic tangency nearby with respect to the horizontal --- or indeed practically any given direction. Along $O^+$, on the other hand, the turns of the isochrons are more gradual and the local maxima due to intersections of the line segment of perturbations are not associated with strong phase sensitivity. Notice further that the penultimate intersection between $\gamma_{0.56} + A \, \mathbf{d}$ and $O^+$ does not come with a nearby quadratic tangency and, hence, does not lead to a local maximum of $\vartheta_{\rm new}$.

%
\section{Phase resetting in a 7D sinoatrial node model}
\label{sec:TKM}
We now illustrate how our computational approach can be applied to systems of dimension higher than two. Indeed, while the multi-segment BVP~\eqref{eq:po}--\eqref{eq:reset} now consists of higher-dimensional subsystems that represent the various orbit segments in this higher-dimensional phase space, the necessary input-output information is still given by the two parameters $\vartheta$ and $\nu$ that determine the relationship $\vartheta_{\rm new} = P(\vartheta_{\rm old})$.

We compute the PTC for the seven-dimensional model from~\cite{TKM} of a sinoatrial node of a rabbit, which is a type of cardiac pacemaker cell. The model is described in standard Hodgkin--Huxley formalism: the main variable is voltage $V$ (measured in ${\rm mV}$), which depends on five ionic currents that are determined by the dynamic opening and closing of six so-called gating variables, denoted $m$, $h$, $d$, $f$, $p$, and $q$. The five currents (measured in ${\rm pA}$) are: a fast inward sodium current $I_{\rm Na}$, a slow inward current $I_{\rm s}$, a delayed rectifier potassium current $I_{\rm K}$, a pacemaker current $I_{\rm h}$, and time-independent leak current $I_{\rm l}$. Then the system of seven equations is given by
\begin{equation}
\label{eq:TKM}
  \left\{ \begin{array}{rcl}
      \dot{V} &=& {\displaystyle - \tfrac{1}{C_m} \left[ I_{Na}(V, m, h) + I_s(V, d, f) + I_K(V, p) + I_h(V, q) + I_l(V) \right]}, \\[2mm]
      \dot{m} &=& {\displaystyle \alpha_m(V) \, (1 - m) - \beta_m(V) \, m}, \\[2mm]
      \dot{h} &=& {\displaystyle \alpha_h(V) \, (1 - h) - \beta_h(V) \, h}, \\[2mm]
      \dot{d} &=& {\displaystyle \alpha_d(V) \,  (1 - d) - \beta_d(V) \, d}, \\[2mm]
      \dot{f} &=& {\displaystyle \alpha_f(V) \,  (1 - f) - \beta_f(V) \, f}, \\[2mm]
      \dot{p} &=& {\displaystyle \alpha_p(V) \, (1 - p) - \beta_p(V) \, p}, \\[2mm]
      \dot{q} &=& {\displaystyle \alpha_q(V) \, (1 - q) - \beta_q(V) \, q}, 
  \end{array} \right.
\end{equation}
where, $C_m = 0.065 \, \mu {\rm F}$ is the capacitance. (Note the minus sign in the right-hand side of the equation for $V$, which was accidentally omitted in~\cite{TKM}.)  The precise form of the ionic currents and the various functions $\alpha_{\rm x}$ and $\beta_{\rm x}$ with ${\rm x} \in \{ m, h, d, f, p, q\}$, and associated parameter values, are given in the Appendix; see also~\cite{TKM}. 

%
\begin{figure}[t!]
  \centering
  \centering 
  \includegraphics{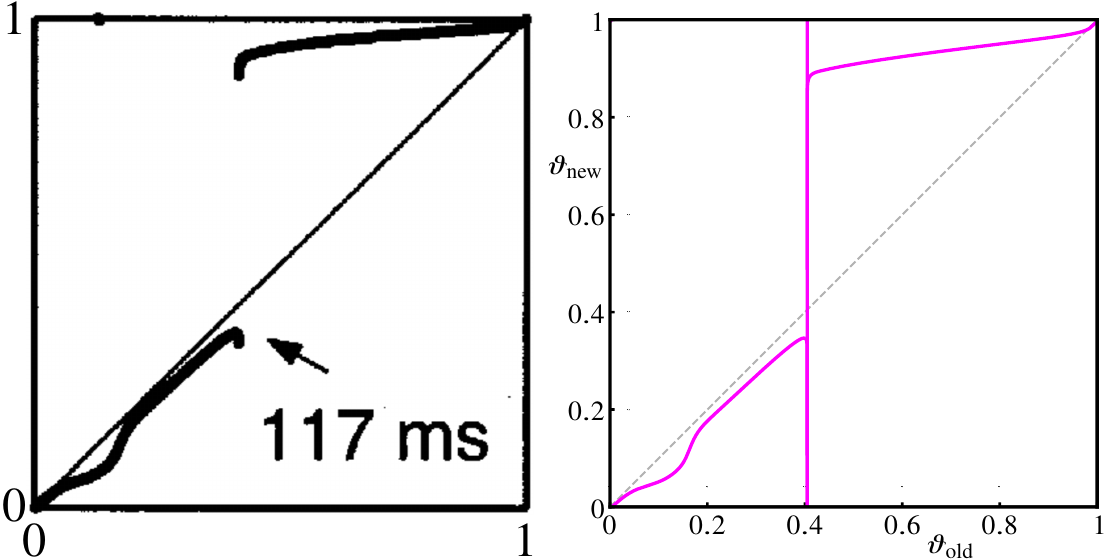}
  \caption{\label{fig:TKM}
    The PTC of the seven-dimensional model~\eqref{eq:TKM}, as presented in~\cite[Fig.~7 (middle)]{TKM}~(a) and as computed with our method~(b). Panel~(a) is from [Krogh-Madsen, Glass, Doedel and Guevara, Apparent discontinuities in the phase-resetting response of cardiac pacemakers, J. Theor. Biol. \textbf{230}(4), 499--519 (2004)] \copyright Elsevier; reproduced with permission.}
\end{figure}
%
%
System~\eqref{eq:TKM} was presented and studied in~\cite{TKM}, because experimental data on similar pacemaker cells suggested that the PTC was discontinuous; see already \fref{fig:TKM}(a). Without a possibility to compute the PTC directly, the authors of~\cite{TKM} reduced the model to a three-dimensional system and used geometric arguments to explain that the apparent discontinuities were abrupt transitions mediated by the stable manifold of a weakly unstable manifold in the model. \Fref{fig:TKM} shows the relevant PTC image from \cite{TKM} and the PTC as computed with our method. The comparison confirms that we are able to calculate the PTC directly in the seven-dimensional model as a continuous curve on $\mT^2$, even though the PTC has a near-vertical segment at $\vartheta_{\rm old} \approx 0.4$.

%
\begin{figure}[t!]
  \centering
  \begin{picture}(150,87)(0,0)
    \put(0,0){\includegraphics{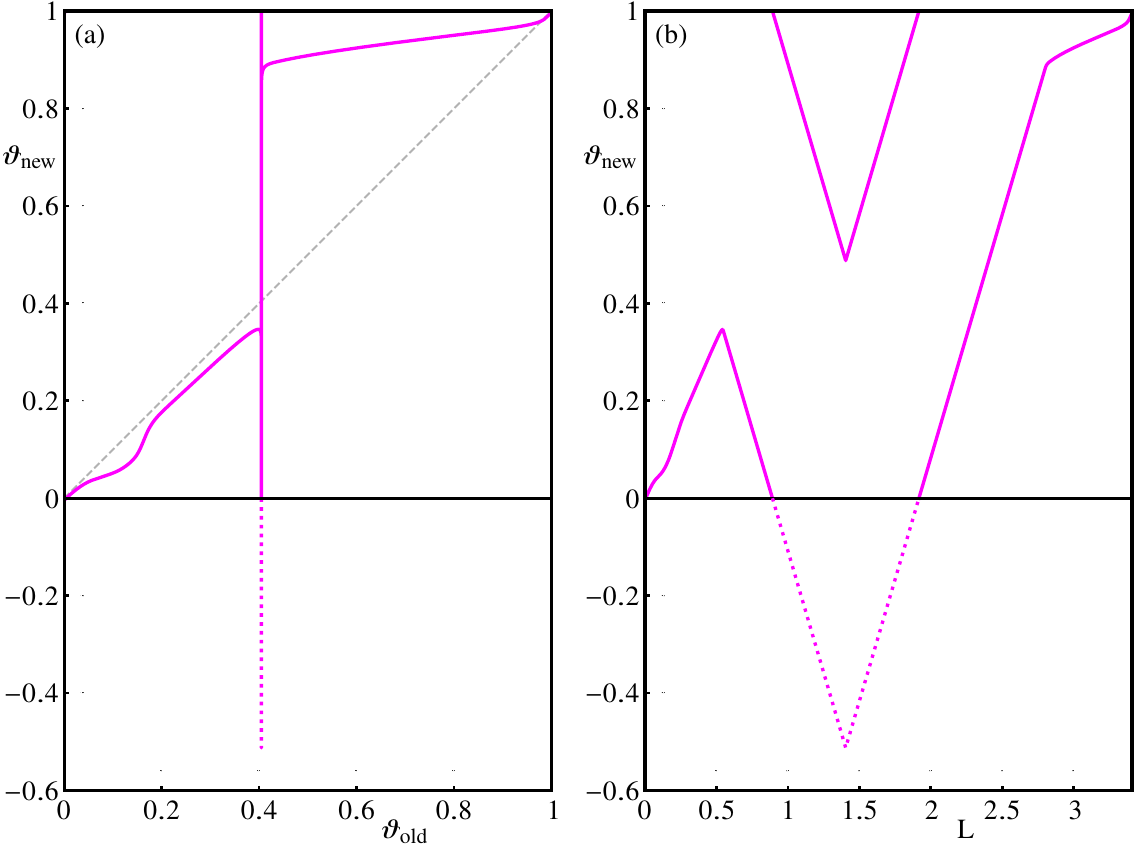}}
  \end{picture}
  \caption{\label{fig:TKMcomp}
    The computed PTC of the seven-dimensional system~\eqref{eq:TKM}. Shown is $\vartheta_{\rm new}$ also over the interval $[-0.6, 0]$ (dotted line) as a function of $\vartheta_{\rm old}$ in panel~(a), and as a function of the arclength $L$ of the PTC in panel~(b).}
\end{figure}
%
%
\Fref{fig:TKMcomp} illustrates that the PTC for system~\eqref{eq:TKM} is indeed continuous. Panel~(a) reproduces the PTC from \fref{fig:TKM}(b), but shows the computed values for $\vartheta_{\rm new}$ over the wider range $[-0.6, 1]$ to show that a maximum of $\vartheta_{\rm new}$ is quickly followed by a minimum of $\vartheta_{\rm new}$ (lowest point of dashed curve). Since it is hard to see that the PTC is indeed continuous, panel~(b) shows $\vartheta_{\rm new}$ over the same range $[-0.6, 1]$, but now as a function of the arclength $L$ of the PTC in the $(\vartheta_{\rm old}, \vartheta_{\rm new})$-plane from the point $(0,\, -3.56 \times 10^{-3})$. The near-vertical segment in the $(\vartheta_{\rm old}, \vartheta_{\rm new})$-plane of panel~(a) correspond to the two (almost) linear segments in the $(L, \vartheta_{\rm new})$-plane of panel~(b). Hence, this representation resolves the steep parts of the PTC in a tiny $\vartheta_{\rm old}$-interval near $0.4$. Panel~(b) also demonstrates that the PTC is indeed a continuous curve closed curve on $\mT^2$ with exactly one maximum at $\vartheta_{\rm new} \approx 0.35$, followed by one minimum at $\vartheta_{\rm new} \approx 0.49$. 

Instead of an instantaneous reset, the reset in~\cite{TKM} is obtained by applying a current with amplitude $I_{\rm app}$ for a fixed duration $\Delta t$; the specific case for which a seemingly discontinuous PTC was observed is given by $I_{\rm app} = -150 \, {\rm pA}$ and $\Delta t = 0.02 \, {\rm s}$. Mathematically, this amounts to replacing the $V$-equation in system~\eqref{eq:TKM} by
\begin{equation}
  \label{eq:Vreset}
  \dot{V} = - \tfrac{1}{C_m} \left[ I_{Na} + I_s + I_K + I_h + I_l \right] + \tfrac{150}{C_m},
\end{equation}
and switching back to the original equation after $\Delta t = 0.02 \, {\rm s}$. In our set-up, this means that we add the perturbation $A \, \mathbf{d}$ to the right-hand side of system~\eqref{eq:TKM}, where the direction vector $\mathbf{d} = (1, 0, 0, 0, 0, 0, 0)$ is the unit vector pointing purely in the $V$-direction and the amplitude $A = 150 / C_m = 2.31$ (${\rm mV} / {\rm s}$).

We include this time-varying perturbation in the multi-segment BVP~\eqref{eq:po}--\eqref{eq:reset} in much the same way as done in~\cite{nota-nc2013}, that is, we replace subsystem~\eqref{eq:svf} defining the orbit segment $\bu$, with boundary conditions~\eqref{eq:scon1E} and~\eqref{eq:scon2}, by two subsystems that define orbit segments $\bu_{\rm ON}$ and $\bu_{\rm OFF}$. Here, $\bu_{\rm ON}$ exists while the applied current is `on' and $\bu_{\rm ON}(1)$ determines the location of the reset~\eqref{eq:reset} after the first $\Delta t = 0.02 \, {\rm s}$. Hence, $\bu_{\rm ON}$ is a solution to system~\eqref{eq:TKM} with equation~\eqref{eq:Vreset} for $V$ with total integration time $\Delta t = 0.02 \, {\rm s}$, that is,
\begin{displaymath}
  \dot{\bu}_{\rm ON} = \Delta t \, \left[ \mathbf{F}(\bu_{\rm ON}) + A \, \mathbf{d} \right].
\end{displaymath}
The second orbit segment $\bu_{\rm OFF}$ is a solution to the original system~\eqref{eq:TKM}, with applied current `off'. The total integration time over both orbit segments combined should be an integer multiple of the period $T_\Gamma$ of the periodic orbit (as before for $\bu$). Hence, we define
\begin{displaymath}
  \dot{\bu}_{\rm OFF} = (k \, T_\Gamma - 0.02) \, \mathbf{F}(\bu_{\rm OFF}).
\end{displaymath}
The subsystem for $\bu_{\rm ON}$ can be viewed as an initial value problem, with initial condition
\begin{displaymath}
  \bu_{\rm ON}(0) = \mathbf{g_u}(0).
\end{displaymath}
Similarly, the initial point of $\bu_{\rm OFF}$ should start where $\bu_{\rm ON}$ ends, that is,
\begin{displaymath}
  \bu_{\rm ON}(1) = \bu_{\rm OFF}(0).
\end{displaymath}
We refer to~\cite{nota-nc2013} for further details.

The end point $\bu_{\rm OFF}(1)$ of the second segment $\bu_{\rm OFF}$ plays the same role as $\bu(1)$ in the multi-segment BVP~\eqref{eq:po}--\eqref{eq:reset}. Hence, $\bu_{\rm OFF}(1)$ must satisfy boundary conditions~\eqref{eq:scon1} and~\eqref{eq:scon2}. Unfortunately, this formulation requires knowledge of the Floquet bundle $\mathbf{v_g}$ specified by subsystem~\eqref{eq:FVE}--\eqref{eq:FVEbc}, and specifically the vector $\mathbf{v_g}(0)$ to measure the distance of $\bu_{\rm OFF}(1)$ to $\Gamma$ in boundary condition~\eqref{eq:scon1}. In the seven-dimensional phase space, this Floquet bundle is no longer unique, because $\Gamma$ now has six non-trivial Floquet exponents. Note that the perpendicular vectors $\bv_0^\perp$, used in boundary conditions~\eqref{eq:atG0} and~\eqref{eq:wuatG0}, and $\mathbf{v_g}(0)^{\perp}$, used in boundary condition~\eqref{eq:scon2}, are still well defined in a higher-dimensional phase space, because the isochrons are codimension-one manifolds. We get around the issue of non-uniqueness as follows. Firstly, we define subsystem~\eqref{eq:FVE}--\eqref{eq:FVEbc} in terms of the adjoint Floquet bundle $\mathbf{v_g}^\perp$, that is, the left eigenvector bundle associated with the trivial Floquet exponent $0$. In other words, we solve the first variational equation
\begin{displaymath}
  \dot{\mathbf{v}}_{\mathbf{g}}^\perp = T_\Gamma \, {\rm D}_{\mathbf{g}} \mathbf{F}^\ast(\mathbf{g}) \, \mathbf{v_g}^\perp,
\end{displaymath}
with the same boundary conditions~\eqref{eq:FVEper} and~\eqref{eq:FVEbc} for $\mathbf{v_g}^\perp$ instead, namely,
\begin{displaymath}
  \left\{ \begin{array}{rcl}
            \mathbf{v_g}^\perp(1) - \mathbf{v_g}^\perp(0) &=& 0, \\
            \mid\!\mid\! \mathbf{v_g}^\perp(0) \!\mid\!\mid &=& 1. 
          \end{array} \right.
\end{displaymath}
Here ${\rm D}_{\mathbf{g}} \mathbf{F}^\ast(\mathbf{g})$ is the transpose Jacobian matrix evaluated along the periodic orbit $\mathbf{g}$. We similarly assume that $\bv_0^\perp$, rather than $\bv_0$, is stored as a known vector.

Secondly, we use the Euclidean norm to measure the distance of $\bu_{\rm OFF}$ from $\mathbf{g}(0)$, that is, we stipulate
\begin{equation}
  \label{eq:scon1E}
  \left[ \bu(1) - \mathbf{g}(0) \right] \cdot \left[ \bu(1) - \mathbf{g}(0) \right] = \eta^2,
\end{equation} 
rather than imposing a signed distance. We remark that a formulation in terms of the Euclidean norm does make the continuation numerically less stable, but it still works for our set-up because $\eta$ is a free parameter that remains positive, and boundary condition \eqref{eq:scon1E} effectively plays a monitoring role.

%
\section{Conclusions}
\label{sec:conclusions}
We presented an algorithm for the computation of the phase reset for a dynamical system with periodic orbit $\Gamma$ that is subjected to an (instantaneous or time-varying) perturbation $\Gamma_A := \Gamma + A \, \mathbf{d}$ of a given direction $\mathbf{d}$ and amplitude $A$. It is well known that the phase reset can be determined from the isochron foliation of the basin ${\cal B}(\Gamma)$, and for small enough $A$, it suffices to know only the linear approximation of the isochrons. Our algorithm tracks the respective nonlinear isochrons and is particularly suited to the computation of phase resets for relatively large $A$. 

Our method is formulated in terms of a multi-segment boundary value problem that is solved by continuation and gives the new phase $\vartheta_{\rm new}$ as a function of either the perturbation amplitude $A$ or the original phase $\vartheta_{\rm old}$ before the reset. The data can readily be used to produce phase transition and phase response curves. We presented the multi-segment BVP set-up in detail for a planar system, but also discussed in \sref{sec:TKM} the straightforward adaptation to higher-dimensional systems, and how to implement phase resets arising from time-varying inputs.

Our approach has the advantage that the map $\vartheta_{\rm old} \mapsto \vartheta_{\rm new}$ is computed in a single continuation run, even in the presence of extreme phase sensitivity. If the amplitude $A$ is such that $\Gamma_A \subset {\cal B}(\Gamma)$, then the associated circle map $P_A : [0, 1) \to [0, 1)$ is obtained in its entirety, and its graph, the phase transition curve (PTC), is a continuous closed curve on $\mathbb{T}^2$. For $A$ close to $0$, the circle map $P_A$ is a near-identity transformation, so that the PTC is a 1:1 torus knot. For large $A$ it is possible that the PTC is a contractible closed curve on the torus, which corresponds to loss of surjectivity of $P_A$. It is well known that surjectivity is lost as soon as $A$ increases past a value for which $\Gamma_A \not\subset {\cal B}(\Gamma)$~\cite{glass84, Winfree2001}. 

There typically exists a maximal amplitude $A_{\rm max}$ such that $\Gamma_A \subset {\cal B}(\Gamma)$ for $0 \leq A \leq A_{\rm max}$. Then $P_A$ depends smoothly on $A$ and, hence, the PTC is a 1:1 torus knot for all $0 \leq A \leq A_{\rm max}$. Therefore, $P_A$ is surjective for all $0 \leq A \leq A_{\rm max}$. We were particularly interested in loss of injectivity of $P_A$ as $A$ increases from $0$. We showed that this is typically mediated by a cubic tangency between the PTC and one of the isochrons of $\Gamma$. Further tangencies lead to very complicated PTCs, with possibly many local maxima and minima and very sudden phase changes. The associated phase sensitivity is known to occur near the boundary of ${\cal B}(\Gamma)$, but our examples illustrate that milder forms of phase sensitivity inside the basin also lead to complicated PTCs. 

We remark that $P_A$ is no longer well defined for all $\vartheta_{\rm old} \in [0, 1)$ when $\Gamma_A$ intersects the boundary of the basin ${\mathcal B}(\Gamma)$. For example, when $\Gamma_A$ crosses an equilibrium that forms a single component of the basin boundary in a planar system, there exists exactly one $\vartheta \in [0, 1)$ such that $P_A(\vartheta)$ is not defined, because the perturbed phase point never returns to $\Gamma$. Entire intervals of $\vartheta_{\rm old} \in [0, 1)$ must be excluded, e.g., when $\Gamma_A$ crosses a repelling periodic orbit of a planar system, such that a closed segment of $\Gamma_A$ lies outside ${\cal B}(\Gamma)$. Changes of the PTC during the transition through different types of boundaries of ${\mathcal B}(\Gamma)$ are beyond the scope of this chapter and will be reported elsewhere.

Phase resets for higher-dimensional systems are expected to exhibit other, more complicated behaviours that lead to possibly different mechanisms of loss of injectivity and/or surjectivity of the circle map associated with the PTC. In particular, the basin ${\cal B}(\Gamma)$ can be a lot more complicated, which affects the isochron foliation and, consequently, the PTC~\cite{Mauroy2015}. Such higher-dimensional systems are of particular  interest when resets are considered in large coupled systems. Even when the coupling is through the mean-field dynamics, such systems can exhibit rich collective dynamics that are reflected in their PTCs~\cite{Bogacz2019, UllnerPoliti2016}. We believe that our approach will be useful in this context, in particular, when the perturbation is a time-dependent stimulus.

%
\section*{Acknowledgments} 
This research is supported, in part by the Royal Society Te Ap\={a}rangi Marsden Fund grant \#16-UOA-286.
We thank Leon Glass and Yannis Kevrekidis for their continued interest in planar isochron computations based on our BVP continuation set-up. Their tiredless enquiries have led to the results presented in this chapter. We thank Michael Dellnitz for stimulating a friendly competitive environment that encouraged us to develop and apply computational methods for invariant manifolds in new contexts.

%
\appendix
\section*{Appendix: Details of the sinoatrial node model} 
System~\eqref{eq:TKM} is the seven-dimensional fast-upstroke model from~\cite{TKM}. The five currents in the equation for $V$ are defined as follows:
\begin{eqnarray*}
  I_{\rm Na} &=& I_{\rm Na}(V, m, h) = g_{\rm Na} \, m^3 \, h \, [V - 40.0], \\[1mm]
  I_{\rm s} &=& I_{\rm s}(V, d, f) = g_{\rm s} \, d \, f \, \left[ e^{(V - 40.0)/25.0} - 1.0 \right], \\[1mm]
  I_{\rm K} &=& I_{\rm K}(V, p) = g_{\rm K} \, p \, \left[ e^{(V + 90.0)/36.1} - 1.0 \right] \, e^{-(V + 40.0)/36.1}, \\[1mm]
  I_{\rm h} &=& I_{\rm h}(V, q) = g_{\rm h} \, q \, \left[ V + 25.0 \right], \\[1mm]
  I_{\rm l} &=& I_{\rm l}(V) \\
         &=& g_{\rm l} \left( 1.2 \left[ 1.0 - e^{ -(V + 60.0)/25.0} \right] + 0.15 \left[ V - 2.0 \right] \, \left[1.0 - e^{-(V - 2.0)/5.0} \right]^{-1} \right).
\end{eqnarray*}
The parameters $g_{\rm Na}$, $g_{\rm s}$, $g_{\rm K}$, $g_{\rm h}$, and $g_{\rm l}$, are conductances (measured in ${\rm nS}$). The capacitance $C_m$ and these five conductances are set to the same values as those of the fast-upstroke model in~\cite{TKM}; see also Table~\ref{tab:TKMpars}.
\begin{table}[t!]
  \begin{center}
    \begin{tabular}{rcllrcllrcll}
      $C_m$ &=& $6.5 \times 10^{-2}$ & ($\mu {\rm F}),$ & \qquad
      $g_{\rm s}$ &=& $1950$ & (${\rm nS}$), & \qquad 
      $g_{\rm h}$ &=& $52$ & (${\rm nS}$), \\
      $g_{\rm Na}$ &=& $325$ & (${\rm nS}$), & \qquad
      $g_{\rm K}$ &=& $354.9$ & (${\rm nS}$), & \qquad 
      $g_{\rm l}$ &=& $65$ & (${\rm nS}$). \\
    \end{tabular}
    \end{center}
    \caption{\label{tab:TKMpars}
      Parameter values for system~\eqref{eq:TKM} as used for the seven-dimensional fast-upstroke model in~\cite{TKM}.}
\end{table}
The $V$-dependent functions for $m$ are defined as 
\begin{displaymath}
  \left\{ \begin{array}{rcl}
    \alpha_m(V) &=& 10^3 \, \left[ V + 37.0 \right] \, \left[ 1.0 - e^{-(V + 37.0)/10.0} \right]^{-1}, \\[3mm]
    \beta_m(V) &=& 4.0 \times 10^4 \, e^{-(V + 62.0)/17.9}, \\[3mm]
    \end{array} \right. 
\end{displaymath}
for $h$, they are defined as
\begin{displaymath}
  \left\{ \begin{array}{rcl}
    \alpha_h(V) &=& 0.1209 \, e^{-(V + 30.0)/6.534}, \\[3mm]
    \beta_h(V) &=& 10^2 \, \left[ e^{-(V + 40.0)/10.0} + 0.1 \right]^{-1}, \\[3mm]
    \end{array} \right. 
\end{displaymath}
for $p$, they are
\begin{displaymath}
  \left\{ \begin{array}{rcl}
    \alpha_p(V) &=& 8.0 \, \left[ 1.0 + e^{-(V + 4.0)/13.0} \right]^{-1}, \\[3mm]
    \beta_p(V) &=& 0.17 \, \left[ V + 40.0 \right] \, \left[ e^{(V + 40.0)/13.3} - 1.0 \right]^{-1}, \\[3mm]
    \end{array} \right. 
\end{displaymath}
for $d$, they are
\begin{displaymath}
  \left\{ \begin{array}{rcl}
    \alpha_d(V) &=& 1.2 \times 10^3 \, \left[ 1.0 + e^{-V/12.0} \right]^{-1}, \\[3mm]
    \beta_d(V) &=& 2.5 \times 10^2 \, \left[ 1.0 + e^{(V + 30.0)/8.0} \right] ^{-1}, \\[3mm]
    \end{array} \right. 
\end{displaymath}
for $f$, they are
\begin{displaymath}
  \left\{ \begin{array}{rcl}
    \alpha_f(V) &=& 0.7 \, \left[ V + 45.0 \right] \, \left[ e^{(V + 45.0)/9.5} - 1.0 \right]^{-1}, \\[3mm]
    \beta_f(V) &=& 36.0 \, \left[ 1.0 + e^{-(V + 21.0)/9.5} \right]^{-1},
    \end{array} \right. 
\end{displaymath}
and finally, for $q$, they are defined as
\begin{displaymath}
  \left\{ \begin{array}{rcl}
    \alpha_q(V) &=& 0.34 \, \left[ V + 100.0 \right] \, \left[ e^{(V + 100.0)/4.4} - 1.0 \right]^{-1} + 0.0495, \\[3mm]
    \beta_q(V) &=& 0.5 \, \left[ V + 40.0 \right] \, \left[ 1.0 - e^{-(V + 40.0)/6.0} \right]^{-1} + 0.0845.
    \end{array} \right.
\end{displaymath}

We now briefly explain how we implement a perturbation, as in the experimental and computed resetting data in~\cite{TKM}, which is obtained by applying a current $I_{\rm app}$ for a fixed duration $\Delta t$. Mathematically, this amounts to replacing the $V$-equation in system~\eqref{eq:TKM} by 
\begin{equation}
  \label{eq:Vreset}
  \dot{V} = {\displaystyle - \tfrac{1}{C_m} \left[ I_{Na} + I_s + I_K + I_h + I_l \right] + 150 \, C_m},
\end{equation}
and switching back to the original equation after $\Delta t = 20 \, {\rm ms}$. In our set-up, this means that we add the perturbation $A \, \mathbf{d}$ to the right-hand side of system~\eqref{eq:TKM}, where the direction vector $\mathbf{d} =(1,0,0,0,0,0,0,0)^T$ is the unit vector pointing purely in the $V$-direction and amplitude $A = 150 \, C_m = 9.75$ (${\rm mV}$). In contrast to the case of instantaneous reset in \sref{sec:2PBVP}, we now replace subsystem~\eqref{eq:svf} defining the orbit segment $\bu$, with boundary conditions~\eqref{eq:scon1E} and~\eqref{eq:scon2}, by subsystems that define two orbit segments $\bu_1$ and $\bu_2$. Here the first orbit segement determines the location of the reset~\eqref{eq:reset} after the first $\Delta t = 20 \, {\rm ms}$. Hence, $\bu_1$ is a solution to system~\eqref{eq:TKM} with equation~\eqref{eq:Vreset} for $V$ with total integration time $\Delta t = 20 \, {\rm ms}$, that is,
\begin{displaymath}
  \dot{\bu}_1 = \Delta t \, \left[ \mathbf{F}(\bu_1) + A \, \mathbf{d} \right].
\end{displaymath}
The second orbit segment $\bu_2$ is a solution to the original system~\eqref{eq:TKM} with an integration time such that the total time over both orbit segments is an integer multiple of the period $T_\Gamma$ of the periodic orbit (as before for $\bu$). Hence, we define
\begin{displaymath}
  \dot{\bu}_2 = (k \, T_\Gamma - 20) \, \mathbf{F}(\bu_2).
\end{displaymath}
The subsystem for $\bu_1$ can be viewed as an initial value problem, with initial condition
\begin{displaymath}
  \bu_1(0) = \mathbf{g_u}(0).
\end{displaymath}
Similarly, the initial point of $\bu_2$ should start where $\bu_1$ ends, that is,
\begin{displaymath}
  \bu_1(1) = \bu_2(0).
\end{displaymath}
Finally, the new phase $\vartheta_{\rm new}$ is selected such that the head point $\mathbf{g_u}(0)$ of the rotated orbit $\mathbf{g_u}$ matches the necessary location imposed by the boundary conditions~\eqref{eq:scon1E} and~\eqref{eq:scon2}  for the end point $\bu_2(1)$ of the second segment $\bu_2$.

%
\bibliographystyle{spphys}
\bibliography{myBibliography}

\end{document}